\documentclass{article}
\usepackage{graphicx}
\usepackage{amsmath, amsthm, amssymb, amscd}
\usepackage{cite}

\def\beq{\begin{equation}}
\def\eeq{\end{equation}}

\def\odd#1{{{d\, }\over{d #1}}}  
\def\ood#1#2{{{d^2 #1}\over{d #2{}^2}}} 
\def\od#1#2{{{d #1}\over{d #2}}} 

\def\pd#1#2{{{\partial #1}\over{\partial #2}}} 

\begin{document}

\title{Geodesics on the Torus\\ and other Surfaces of Revolution\\ Clarified Using Undergraduate Physics Tricks\\
 with Bonus: Nonrelativistic and Relativistic Kepler Problems
}
\author{
Robert T. Jantzen\\ 
Department of Mathematics and Statistics\\
Villanova University
}
\date{May 23, 2010}
\maketitle

\begin{abstract}
In considering the mathematical problem of describing the geodesics on a torus or any other surface of revolution, there is a tremendous advantage in conceptual understanding that derives from taking the point of view of a physicist by interpreting parametrized geodesics as the paths traced out in time by the motion of a point in the surface, identifying the parameter with the time. Considering energy levels in an effective potential for the reduced motion then proves to be an extremely useful tool in studying the behavior and properties of the geodesics. The same approach can be easily tweaked to extend to both the nonrelativistic and relativistic Kepler problems. The spectrum of closed geodesics on the torus is analogous to the quantization of energy levels in models of atoms.
\end{abstract}

\section{Geometry as Forced Motion?}

The central idea of general relativity \cite{wikiGR,walecka} is that motion under the influence of the gravitational force is viewed instead as the natural result of geodesic motion in a curved spacetime. It was this application to physics that gave a tremendous boost to interest in classical differential geometry at the turn of the twentieth century just as it was maturing as a mathematical discipline. The term \lq geodesic' is not part of our common everyday language, but it is a simple idea that most of us have some intuition about from the role of straight lines in flat Euclidean geometry and great circles on a sphere: geodesics are paths which at least locally minimize the distance between two points in a space and are such that the direction of the path within the space does not change as we move along it. For example, moving along a great circle on a sphere, the direction of the path only changes within the enveloping 3-dimensional Euclidean space to stay tangent to the sphere, but within the sphere itself, does not rotate to the left or to the right with respect to its intrinsic geometry. If one imagines an insect-sized toy car on a much larger sphere, moving along a great circle means locking the steering wheel so that the car always moves straight ahead, never changing direction.

The torus (donut shaped surface) is a nice example to follow the flat plane and the highly symmetric curved sphere as a playground for gaining concrete experience with geodesics \cite{irons}, but the techniques discussed here which prove useful in its description apply to any surface of revolution \cite{carmo,gray3,oprea,pressley}, provided that the curve being revolved has an explicit arclength parametrization (not true of a simple parabola of revolution, for example). These techniques are familiar to any physics undergraduate student who has taken a second course in classical mechanics, the motivating application which for the most part led to the development of calculus itself. These elementary tools are not part of the typical mathematician's toolbox, however, so the description of this problem one finds in the literature or on the web falls short of what it should be. 

One imagines tracing out a parametrized geodesic path in a surface as the motion of a point particle in the surface, interpreting the parameter in terms of which the curve is expressed as the time. Indeed a computer algebra system animation of a parametrized curve does exactly this in some appropriately chosen time unit. We will limit ourselves to consider surfaces of revolution about the $z$-axis in ordinary Euclidean space, so that one can describe the surface using an ``azimuthal'' angular variable $\theta$ measuring the angle about the $z$-axis (just the usual polar coordinate of the projection of a point into the $x$-$y$ plane), and a radial arc length variable $r$ which parametrizes the curve of constant angle $\theta$ in the surface by its arc length (not to be confused with the usual polar coordinate of the projection into that plane).
Like the Kepler problem for orbits in a plane expressed in polar coordinates, or more generally for motion in any central force field in 2-dimensions, conservation of angular momentum about the axis of symmetry enables one to reduce the problem to that of the radial motion alone with the effects of the angular motion felt by an effective potential, the centrifugal potential, which acts as a barrier surrounding the axis of symmetry for motion with nonzero angular momentum. One glance at the diagram of this 1-dimensional effective potential captures the key features of the geodesic problem. Yet mathematicians appear to be unaware of this possibility. Visualization of mathematics is one of the most helpful tools for understanding, so it is worth bringing attention to it in this context. 

With simple modifications one can also incorporate a central force into this geodesic problem in 2 dimensions and apply the same machinery to reveal the conic section orbits of the Kepler problem, i.e., the nonrelativistic gravitational problem, as well as understand its relativistic generalization to motion of a test particle in a nonrotating black hole or outside any static spherically symmetric body. One does not even really need the advanced classical mechanical knowledge (the calculus of variations and Lagrangians, see Appendix B) as long as one buys into the concept of the conservation of energy and of angular momentum, ideas already presented in high school physics.

In this way we can reintroduce our even limited intuition about particle motion into the problem of finding geodesics on a surface of revolution, thus converting geometry back into the framework of forced motion. For an undergraduate who has taken separately multivariable calculus, differential equations and linear algebra, a little exposure to elementary differential geometry brings all of these tools together under one roof and has the potential to play a unifying role in integrating them together, something sorely lacking in many undergraduate mathematics programs for students who want or need some working knowledge of those topics. Using geodesic motion on the torus as a central example gives an elegant and concrete focus to this particular class of problems. As a classical mechanical system in physics, it is a particularly rich and beautiful example of motion in a central force field in 2 dimensions, a consequence of its rotational and discrete symmetries.

In the present article, we suggestively introduce the main ideas needed to consider geodesics, but this cannot substitute for a more detailed study of the necessary mathematical tools. Thus it is really aimed at those who already have enough background, but does not exclude those with interest in the topic who are willing to accept the loose explanation given of the preliminary ideas. 

Finally, searching the literature one finds very little concrete discussion of the geodesics on the torus precisely because they cannot be fully described analytically, but now that computer algebra systems are readily available, one begins to see some limited mention of using computer algebra systems (Maple and Mathematica) to study them numerically. However, even if one accepts the defining differential equations as given and plays with numerical solutions, it is still important to have a good understanding of the setting in which these numerical games are played. The present discussion hopes to supply a detailed picture that is still missing elsewhere.
In particular, the classification of closed geodesics on the torus corresponds to a discrete set of energy levels in this picture, mirroring the analogous quantization of energy levels in the model of an atom. Finding and studying these closed geodesics is equally interesting to those of us who enjoy mathematics for its own sake.

\section{Straight lines, circles, spheres and tori}

For straight lines in the $x$-$y$ plane, the arc length of a segment is simply the distance between the points, which amounts to the Pythagorean theorem. To describe the differential of arc length between two nearby points on a curve in the plane, the same Pythagorean theorem is used
\beq
  ds^2 = dx^2 + dy^2 \,.
\eeq
For polar coordinates in the plane
\beq
  x = r \cos\theta\,,\
  y = r \sin\theta\,,\
\eeq
since the radial and angular coordinate lines are orthogonal, one can evaluate the differential of arc length along a curve by applying the Pythagorean theorem to the two orthogonal differentials of arc length along the two coordinate lines. $dr$ directly measures the differential of arc length in the radial direction, while $r \,d\theta$ measures the differential of arc length in the angular direction, since increments of angle must be multiplied by the radius of the arc to obtain the arc length of a circular arc. Thus one obtains
\beq
  ds^2 = dr^2 + r^2 d\theta^2\,,
\eeq
a result which could also be obtained simply by taking the differentials of $x$ and $y$ expressed in terms of $r$ and $\theta$ and substituting into $ds^2$, then expanding and simplifying the result using the fundamental trigonometric identity.

If we repeat this in the context of Cartesian coordinates $(x,y,z)$ in Euclidean space where $ds^2=dx^2+dy^2+dz^2$, the polar coordinates in the plane are upgraded to cylindrical coordinates in space
\beq
   x = r \cos\theta\,,\
  y = r \sin\theta\,,\
  z = z\,,
\eeq
and one passes to spherical coordinates by introducing polar coordinates in the $r$-$z$ half plane $r\ge 0$
\beq
  r = \rho \sin\phi\,,\
  z = \rho \cos\phi\,,\
\eeq
so that taken together one has
\beq
  x = \rho \sin\phi \cos\theta\,,\
  y = \rho \sin\phi \sin\theta\,,\
  z = \rho \cos\phi\,.
\eeq
Since these coordinates are also orthogonal, the iterated Pythagorean theorem (distance formula as a sum of squares) can be used to evaluate the differential of arc length along a curve. $\rho$ is an arc length coordinate, while $r \,d\theta = \rho\sin\phi\, d\theta$ continues to describe the differential of arc length along the azimuthal angular coordinate ($\theta$) lines, and now $\rho\, d\phi$ describes the differential of arc length along the polar angular coordinate ($\phi$) lines, so
\beq
 ds^2 = d\rho^2 + \rho^2 d\phi^2 + \rho^2\sin^2\phi\, d\theta^2\,. 
\eeq
This same result could have been obtained by inserting the differentials of the Cartesian coordinates into $ds^2$, expanding and simplifying.

For a sphere of constant radius $\rho=b$, this reduces to 
\beq
 ds^2 = b^2 d\phi^2 + b^2\sin^2\phi\, d\theta^2
      = [d(b\,\phi)]^2 + [b\sin\phi]^2 d\theta^2
      = dr^2 +  [b\sin(r/b)]^2 d\theta^2 \,, 
\eeq
where $r = b\,\phi$ is a new arc length coordinate measuring the arc length along the $\theta$ coordinate lines down from the North Pole of the sphere, not to be confused with the previous polar coordinate $r$. We do this to compare the sphere arc length formula directly with the plane formula in the same variables.
If we had instead introduced an angle $\varphi=\pi/2-\phi$ in the constant $\theta$ half plane measured up from the horizontal direction so that $\sin\phi=\cos\varphi$ and then a new radial coordinate $\tilde r=b\,\varphi$, we would have found instead
\beq
 ds^2 = d\tilde r^2 +  [b\cos(\tilde r/b)]^2 d\theta^2 \,. 
\eeq
Thus for both the flat plane and the sphere, a common form of the so called \lq\lq line element'' $ds^2$ expressed in terms of a single function $R(r)$ which gives the radius of the azimuthal $\theta$ coordinate circle corresponding to the fixed radial variable $r$ within the surface
\beq
 R(r) = \begin{cases}
        r\,, &  \text{polar coords in plane,}\cr
        b^{-1}\sin(r/b)\,, &  \text{usual spherical coords ($r/b= \phi$),}\cr
        b^{-1}\cos(r/b)\,, &  \text{alternative spherical coords ($r/b= \pi/2-\phi$)}\cr
        \end{cases}
\eeq
and the standard form of the metric
\beq
 ds^2 = dr^2 +  R(r)^2 d\theta^2 \,. 
\eeq

At this point I must remark that spherical coordinates trigger my mathematics/physics schizophrenia. As a physicist I normally use the opposite convention for spherical coordinates \cite{sphericalcoords} common in physics applications in which $(\theta,\phi)$ are switched with respect to the usual calculus textbook convention adopted here, so that $\phi$ is the azimuthal angle in cylindrical and spherical coordinates rather than $\theta$ as in polar coordinates. To make matters worse, the radial coordinates $(r,\rho)$ are also interchanged in the physics convention! As a calculus teacher myself, I have to admit that it is easier to not switch the angles on students who have enough to keep track of as it is, and in the present context it is easiest to treat all of these problems with the same notation. 

To assign unique polar, cylindrical or spherical coordinates to a point in the plane or space, one restricts the radial variables $r$ and $\rho$ to be nonnegative and the zenith angle $\phi$ to its obvious closed interval, but one has the option of choosing one of two obviously useful intervals for the azimuthal angle $\theta$
\begin{eqnarray}
r\ge0:&& \mbox{radial distance from $z$ axis,}\nonumber\\
\rho\ge0:&& \mbox{radial distance from origin,}\nonumber\\
0< \phi \le \pi:&& \mbox{zenith angle down from positive $z$ axis,}\nonumber\\
0\le \theta <2\pi &\mbox{or}&\ -\pi< \phi \le \pi:\qquad
 \mbox{azimuthal angle.}
\end{eqnarray}
However, to describe in a continuous fashion those parametrized curves which wrap more than one complete revolution around the $z$ axis, one must allow the azimuthal angle to take values outside these intervals.

Both the flat plane $z=0$ within space and the spherical coordinate sphere of radius $b$ arise as surfaces of revolution \cite{carmo,gray3,oprea,pressley}, which are constructed by taking a curve in the $x$-$z$ plane and revolving it around the $z$-axis to obtain a surface of revolution. As long as the curve can be parametrized so that the integral defining the arc length parametrization can be evaluated exactly, and one can invert the relation between the arc length and curve parameter to express the curve in terms of an arc length parametrization, we can play the same game as above. Revolving the $x$-axis itself leads to the flat plane $z=0$, while revolving the circle $x^2+z^2=b^2$ leads to the spherical coordinate sphere of radius $b$. Since arc length is trivial for straight lines and arcs of circles, choosing any straight line or circle in the $x$-$z$ plane will do the job. An inclined straight line leads to a flat cone (for example, a $\phi$ coordinate surface in spherical coordinates, where $r=\rho$ and $R(r)= r\sin\phi_0$), while a vertical straight line leads to a flat cylinder  (for example, an $r$ coordinate surface in cylindrical coordinates, where the new radial coordinate is $r=z$ and $R(r)= 1$). Choosing instead a circle not centered on the $z$-axis leads to a torus.

If the curve of this construction crosses the $z$-axis and its tangent line is not horizontal there, then the surface of revolution has a singular point where it is not differentiable, and the limiting tangent line if not vertical either sweeps out a cone with its vertex at this point, or corresponds to a limiting cusp of revolution if indeed it is vertical. For a torus, this conical case leads to self-intersections on the axis of symmetry, with two disjoint components of the surface.

\begin{figure}[t] 
\typeout{*** EPS figure torus1}
\begin{center}
\includegraphics[scale=0.3]{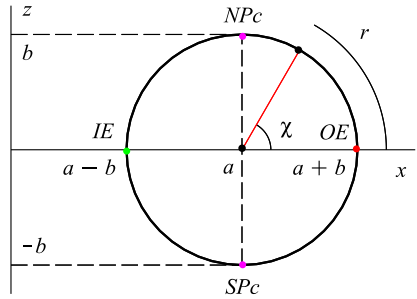}\qquad
\includegraphics[scale=0.4]{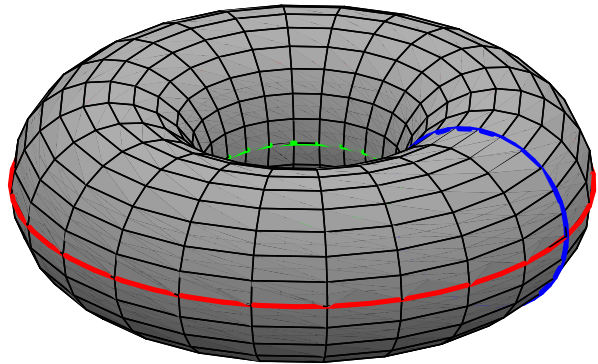}
\end{center}
\vglue-1cm
\caption{
The torus is a surface of revolution obtained by revolving about the $z$-axis a circle in the $x$-$z$ plane with center on the $x$-axis.  Illustrated here is the ``unit ring torus" case $(a,b)=(2,1)$ of a unit circle which is revolved around the axis, with an inner equator of unit radius. The outer ($\chi=0$, red) and inner ($\chi=\pm\pi$, green) equators are shown together with the ``prime meridian" ($\theta=0$, blue). The Northern ($\chi=\pi/2$) and Southern  ($\chi=-\pi/2$) Polar Circles correspond to the North and South Poles on the sphere. The radial arc length coordinate $r=b\,\chi$ and the corresponding angle $\chi$ are measured upwards from the outer equator. The grid shown in the computer rendition of the surface consists of the constant $\theta$ circles which result from the intersection of the torus with vertical planes through the symmetry axis (the meridians) and the constant $r$ or $\chi$ horizontal circles (the parallels) which result from the intersection of the torus with horizontal planes.
} 
\label{fig:torus1}
\end{figure}

We construct a torus as a surface of revolution by revolving a circle about the $z$ axis
as illustrated in Fig.~\ref{fig:torus1}.
We take a circle of radius $b>0$ with center at $x=a>0$ on the $x$-axis, parametrized by a polar angle $\chi$ measured in the $x$-$z$ plane counterclockwise from the positive $x$-axis,
and rotate this circle about the $z$-axis (the ``symmetry axis").
We thus obtain a standard family of tori for which the radius of the azimuthal radius function is 
$ R=a+b\cos\chi$,
while $z=b\sin\chi$ so that the Cartesian coordinates are parametrized by
\beq
  x = (a+b\cos\chi) \cos\theta\,,\
  y = (a+b\cos\chi) \sin\theta\,,\
  z =    b\sin\chi\,,
\eeq
or in vector form, introducing the position vector $\vec r = \langle x, y, z \rangle$
\beq
 \vec r = \langle (a+b\cos\chi) \cos\theta,(a+b\cos\chi) \sin\theta, b\sin\chi 
          \rangle \,.
\eeq
Then introducing the arc length coordinate $r=b\chi$, this becomes 
\beq
 \vec r = \langle (a+b\cos(r/b)) \cos\theta,(a+b\cos(r/b)) \sin\theta, b\sin(r/b) 
          \rangle \,,
\eeq
and we can add to the above list of azimuthal radius functions
\beq
R(r)= a+b\cos(r/b) \,,\quad  \rm torus\,,
\eeq
which reduces to the alternative spherical coordinate relation when $a=0$. As before one could also derive the standard form for the line element by merely substituting into $ds^2$ the differentials of the Cartesian coordinates expressed in terms of the two angular variables, expanding and simplifying and rescaling the radial angular variable to $r$, but since the coordinate circles are orthogonal, it is just the sum of the squares of the simple circular differentials of arc length along the two directions. 

\clearpage

By introducing the shape parameter $c=(a-b)/b\ge-1$ 
in terms of which one has $a/b=1+c$ and $(a+b)/b=2+c$, 
one finds the following distinct subsets of the whole family. The ring tori \cite{torus-online} 
result from the range $a>b$ or $c>0$.
A horn torus results from $a=b$ or $c=0$ in which the inner radius of the torus goes to zero, and the origin itself has a limiting behavior which is the tip of two cusps of revolution lacking a tangent plane so it is a bit problematic.
For $0<a<b$ or $-1<c<0$ one has spindle tori with two self-intersection points on the axis of symmetry and an additional part of the surface inside the outer surface.
For the limiting case $a=0$ or $c=-1$ one has the ordinary sphere of radius $b$ where $\chi=\varphi$.
Note that apart from on overall scale factor which doesn't change any geodesic properties of the torus (which only depend on the shape of the torus), we have a one-parameter family of distinct (nonsimilar) torus geometries parametrized by the shape parameter $-1< c <\infty$,
with $c=-1$ corresponding to the limiting case of the sphere. 
Note that the commonly used horn and spindle terminology is interchanged by Gray et al \cite{gray3} for some reason, while apparently Kepler named the inner and outer surface components of the self-intersecting spindle tori to be lemons and apples respectively for obvious shape reasons \cite{wein}.

Note that if one expresses the differential of arc length in terms of the two angular coordinates and the shape parameter, one finds
\beq
 ds^2 = b^2 [ d\chi^2 + (c+1 +\cos\chi)^2 d\theta^2 ]\,. 
\eeq
The overall scale factor $b^2$ just rescales the unit of length on the torus, but the so called ``conformal geometry" only depends on the shape parameter $c$.

To have a concrete simple ring torus example to play with, it seems natural to use the torus with $(a,b)=(2,1)$ and shape parameter $c=1$, which we will call the ``unit ring torus" (my terminology). Any ring torus has a Northern and Southern Polar Circle (my terminology but obvious) like the North and South Poles on a sphere, and both an inner and outer equator like the equator on a sphere. Like the sphere, its plane cross-sections at constant $\theta$ are circles, ``the meridians," marking off longitude from the ``prime meridian" $\theta=0$, whose intersection with the outer equator is the origin of the $(r,\theta)$ coordinate system. Similarly each circle of revolution in the torus is a ``parallel" or line of latitude as on the sphere. This terminology of meridians (plane cross-sections through the axis of symmetry) and parallels (circles of revolution about the axis) introduced for any surface of revolution also applies to the non-ring tori. The inner equator shrinks to a point at the origin for the horn torus and then passes to an outer equator of the inner lemon surface for the spindle tori. We will refer to the horn torus with $(a,b)=(1,1)$ and $c=0$ as the ``unit horn torus."
The upper/lower half of any torus above/below the equators might be called the upper/lower hemitorus in analogy with the hemispheres of a sphere, but one can also introduce the inner/outer hemitorus inside/outside the Northern and Southern Polar Circles.

One can use as convenient range for either angular coordinate $\theta$ or $\chi$ the interval $[0,2\pi)$ or $(-\pi,\pi]$ but for $r$ we will usually use the range $(-b\,\pi,b\,\pi]$. One can also conveniently map the pair of angles onto the unit square $[0,1]\times [0,1]$ by measuring both angles in units of revolutions (divide each angle in the range $[0,2\pi)$ by $2\pi$), identifying opposite edges of the square. This can be useful when one is considering closed geodesics, where this will tell us information about the horizontal and vertical periods of the geodesics compared to the fundamental period $2\pi$ shared by both angular coordinates. 
Such closed geodesics can then be characterized by a pair of integers $(m,n)$, given a fixed reference point through which they pass: the number $m$ of vertical oscillations during $n$ azimuthal revolutions around the symmetry axis (without retracing their path). Note that while all of the various coordinates are periodic functions of an affine parameter along a geodesic, only for a closed geodesic does a common period exist so that the geodesic itself as a function from the real line to the torus manifold is periodic.

We will see that the inner and outer equators and the meridians are all geodesics, but for the ring tori on which we will concentrate our attention, the only geodesic which does not pass through the outer equator is the inner equator. With this one exception, one can classify all geodesics on a ring torus by studying the initial value problem at a fixed point on the outer equator. The meridian passing through the initial point turns out to be a geodesic. The remaining geodesics through this point have an obvious categorization into two separate families: those which cross the inner equator and thus require an unbounded range of values of the radial coordinate to describe 
(the ``unbound geodesics") and the rest,
for which the range of values of the radial coordinate is bounded
 (the ``bound geodesics"). There are closed geodesics of each type so one needs three integers $[m,n;p]$ to completely characterize them, where $p=0$ refers to the bound case and $p=1$ to the unbound case with one exception. The bound terminology is standard in the physics approach to the problem to be described below.
 The inner and outer equators are closed and correspond respectively to the triplets $[0,1;1]$ and $[0,1;0]$ if we agree to associate the bound inner equator geodesic with all the unbound geodesics which cross it.  Each meridian is a closed unbound geodesic corresponding to $[1,0;1]$. The case $[1,0;0]$ does not exist. For the horn torus, the inner equator degenerates to a point and there are only bound geodesics. For the spindle tori there are instead two disjoint families of bound geodesics which can be labeled by $p=0$ for those confined to the outer apple surface and $p=1$ for those confined to the inner lemon surface. In this latter case one must consider two separate initial value problems to describe all the geodesics, one on each equator.

Although it is pretty obvious, it is worth noting that in addition to its rotational symmetry,
each torus is reflection symmetric about the equatorial plane and every vertical plane through the axis of symmetry. This will also be true of its geodesics. In particular for initial data on the outer equator, these reflection symmetries of the torus correspond to reflections across the radial (vertical) and azimuthal (horizontal) directions in the tangent space, so
it is enough to consider initial angles in one quadrant of the vertical tangent plane to study all types of geodesics which pass through that initial point. Because of the rotational symmetry about the vertical axis, it does not matter where on the outer equator we consider the initial value problem. 

It is important to also note that the continuous symmetry group of this problem combined with its discrete symmetries gives the dynamical system of geodesic motion on the torus a rich structure making possible a great deal of insight that a system without symmetry does not allow. This is ultimately responsible for most of the analysis of the present article, and is usually missing in differential geometry textbooks dealing with the classical topics of curves and surfaces in space. 

\section{Geodesic equations}
 
Let $(x^1,x^2)=(r,\theta)$ denote the two parameters/coordinates on the surface of revolution in general and on the torus in particular. The position vector as a function of these variables is then of the form
\beq
\vec r(r,\theta) 
= \langle x(r,\theta),y(r,\theta),z(r,\theta)\rangle
= R(r) \langle \cos\theta,\sin\theta,0 \rangle + Z(r) \langle 0,0,1 \rangle \,,
\eeq
where $Z(r)$ is the height function with respect to the $x$-$y$ plane, while $R(r)$ is the horizontal distance from the $z$-axis. Considering each of the variables in turn as a parameter for a curve with the other variable fixed leads to the two orthogonal tangent vectors along the $r$ and $\theta$ coordinate lines on the surface
\begin{eqnarray}
  e_r &=& \pd{\vec r(r,\theta)}{r} 
       = R'(r)  \langle \cos\theta,\sin\theta,0 \rangle + Z'(r) \langle 0,0,1 \rangle
\,,\\
  e_\theta &=& \pd{\vec r(r,\theta)}{\theta}
            = R(r) \langle -\sin\theta,\cos\theta,0 \rangle
\,,
\end{eqnarray}
where for the torus one has
\beq
 R'(r)=-\sin(r/b)\,,\
 Z'(r)=\cos(r/b) \,.
\eeq
The  inner products of these two tangent vectors give the so called metric components
$g_{ij} = e_i \cdot e_j$ ($i,j=r,\theta$)
whose values are given by
\beq
  g_{rr} = R'(r)^2 +Z'(r)^2 = 1 \,,\quad
  g_{\theta\theta} = R(r)^2 \,,\quad
  g_{r\theta}=g_{\theta r} = 0 \,.
\eeq
The first result $g_{rr}=1$ is a consequence of the fact that by choice $r$ is an arc length parameter along its coordinate lines. The vanishing components reflect the orthogonality of the two tangent vectors, which form a basis of the 2-dimensional tangent space.
The matrix of metric components with the indices $(r,\theta)$ used interchangeably with $(1,2)$ is diagonal because of this orthogonality.
\beq
        \left(\begin{array}{cc} g_{rr}&g_{r\theta}\\ g_{\theta r} & g_{\theta\theta}\\ \end{array} \right)
  =  \left(\begin{array}{cc} 1&0\\ 0 & R(r)^2\\ \end{array} \right) \,.
\eeq

The differential of the position vector is then
\beq
  d \vec r(r,\theta) = e_r \,dr + e_\theta\, d\theta 
\eeq
and its self-dot product defines the so called line element representing the square of the differential of arc length on the surface
\beq
  ds^2 =  d \vec r(r,\theta) \cdot  d \vec r(r,\theta) \,,
\eeq
analogous to the expressions in polar or spherical coordinates in the flat plane or in flat space.
This line element is a quadratic form in the differentials of these parameters/coordinates whose symmetric matrix of coefficients is referred to as the components of the metric. Using the Einstein summation convention that repeated indices are summed over, and using both $1,2$ and $r,\theta$ for the indices $(i,j)$ as convenient, one therefore has a line element and metric component matrix related by
\begin{eqnarray}
  ds^2 &=& dr^2 +  R(r)^2 d\theta^2 
  = \sum_{i=1}^2 \sum_{j=1}^2 g_{ij} dx^i dx^j
  =  g_{ij} dx^i dx^j
\\
  &=& g_{rr} dr^2
   +g_{r\theta} dr d\theta
   +g_{\theta r} d\theta dr 
   +g_{\theta\theta} d\theta^2 
\,, \nonumber
\end{eqnarray}
where
\beq
    g_{rr} = 1\,,\
  g_{r\theta} =
   g_{\theta r} = 0\,,\
  g_{\theta\theta} = R(r)^{2}
\,.
\eeq

The inverse of the matrix of metric components has components denoted by $g^{ij}$, which for this diagonal matrix is just the diagonal matrix of reciprocals of the original diagonal entries
\beq
    g^{rr} = 1\,,\
  g^{r\theta} =
   g^{\theta r} = 0\,,\
  g^{\theta\theta} = R(r)^{-2}
\,,
\eeq
namely
\beq
        \left(\begin{array}{cc} g^{rr}&g^{r\theta}\\ g^{\theta r} & g^{\theta\theta}\\ \end{array} \right)
  =  \left(\begin{array}{cc} 1&0\\ 0 & R(r)^{-2}\\ \end{array} \right) \,.
\eeq

It is useful to normalize the two tangent vectors along the coordinate lines
\beq
  e_{\hat i} = g_{ii}^{-1/2} e_i :\qquad
  e_{\hat r} = e_r\,,\
  e_{\hat \theta} = R(r)^{-1}e_\theta
\eeq
to obtain an orthonormal basis of the tangent plane to the surface 
\beq
  e_{\hat r} \cdot  e_{\hat r} =1 =  e_{\hat \theta} \cdot  e_{\hat \theta}\,,\quad
  e_{\hat r} \cdot e_{\hat \theta} =0 \,.
\eeq 
Then any vector $X$ tangent to the torus can be expressed in the form
\beq
  X = X^r e_r + X^\theta e_\theta
    = X^{\hat r} e_r + X^{\hat\theta} e_{\hat\theta} \,,
\eeq
where the two components $X^r$ and $X^\theta$ are functions of the two coordinates.
We can also evaluate a unit normal vector field to the surface which is the outward unit normal on the torus by taking the cross-product of these two orthonormal vectors reversed in sign
\beq
   \hat n = -e_{\hat r} \times e_{\hat\theta} = -R(r)^{-1} e_r \times e_\theta
          =  Z'(r) \langle \cos\theta,\sin\theta,0\rangle
            -R'(r) \langle 0,0,1\rangle
\,,
\eeq
where the radial arc length condition again confirms that this is clearly a unit vector.

Since the two basis vectors $e_r$ and $e_\theta$ not only change in directions tangent to the surface as one moves around in the surface, but also must change in the normal direction to remain tangent to the surface as the orientation of the tangent plane changes,
their partial derivatives with respect to $r$ and $\theta$ will consist of a part tangent to the surface and a part along the normal vector to the surface. Thus these derivatives can be expressed as linear combinations of $e_r$ and $e_\theta$ plus a multiple of the unit outward normal $\hat n$. 
\beq
  \frac{\partial e_i}{\partial x^j} = \Gamma^k{}_{ij} e_k + K_{ij} \hat n
\,.
\eeq
The expansion coefficients $\Gamma^k{}_{ij}$ of the part tangent to the surface are called the Christoffel symbols. Ignoring the extra term along the normal direction, appropriate if we are only interested in how vectors change within the surface, leads to the so called covariant derivative within the surface.  This is useful since we want to define a geodesic curve by the property that its tangent vector does not change its direction  within the surface (or its length if properly parametrized) as we move along the curve. The additional term along the unit normal is a consequence of the bending of the surface, and its coefficients are the components of the ``extrinsic curvature" of the surface, also called the second fundamental form modulo a sign convention. Here we are only interested in the intrinsic geometry of the surface, so we will not worry about this object, interesting in its own right.

We can define the covariant derivative within the surface to be the ordinary partial derivative on scalar functions in the surface (functions of $x^i$, i.e., of $(r,\theta)$, functions which can arise by re-expressing functions of the original Cartesian coordinates on space, or which are simply new functions of the surface coordinates)
\beq
\nabla_i f = \frac{\partial f}{\partial x^i}\,,
\eeq 
but for vector fields tangent to the surface we can simply ignore the partial derivative component along the normal direction by defining the covariant derivatives of the two basis vector fields using the Christoffel terms only
\beq
\nabla_i e_j = \Gamma^k{}_{ij} e_k\,.
\eeq
We can then extend them to a linear combination of the two surface vector fields by linearity and the product rule to the following linear combination of products of scalars and basis vector fields
\begin{eqnarray}
  X &=& X^r e_r + X^\theta e_\theta = X^i e_i\,,
\\
  \nabla_j X 
 &=& (\nabla_j{X^i}) e_i 
  + X^i (\nabla_j e_i) 
 = \frac{\partial X^i}{\partial x^j} e_i 
  + X^i\Gamma^k{}_{ji} e_k 
 = \frac{\partial X^i}{\partial x^j} e_i 
  + X^k\Gamma^i{}_{jk} e_i 
\nonumber\\
 &=& \left(\frac{\partial X^i}{\partial x^j} e_i 
  + \Gamma^i{}_{jk} X^k\right) e_i 
 \equiv [\nabla_{e_j} X^i ]e_i \,,
\end{eqnarray}
after a convenient relabeling of the indices which are summed over.
The extra correction terms just take into account that not only are the components of the vector $X$ changing, but also the basis vectors themselves. Furthermore, this only measures changes of the vector $X$ within the surface, ignoring how it must change in the normal direction to remain tangent to the surface.

Since the basis vectors are orthogonal to themselves and $\hat n$, we can easily project out the Christoffel symbol components by taking dot products of their defining relation
\beq
  e_\ell \cdot   \frac{\partial e_i}{\partial x^j} 
       = \Gamma^k{}_{ij} e_\ell \cdot e_k
       = g_{\ell k} \Gamma^k{}_{ij} 
        =\sum^2_{k=1} g_{\ell k} \Gamma^k{}_{ij} 
        = g_{\ell \ell} \Gamma^\ell{}_{ij} \quad \mbox{(no sum on $\ell$)}
\,,
\eeq
where we remind the reader that the repeated index $k$ is summed over, but each sum only consists of one term since $g_{\ell k}=0$ if $k\ne \ell$.
It turns out that the Christoffel symbols are also given directly by derivatives of the metric components
\beq
 \Gamma^i{}_{jk} 
= \frac12 g^{i\ell}\left( \pd{g_{j\ell}}{x^k}+ \pd{g_{k\ell}}{x^j} - \pd{g_{jk}}{x^\ell}
         \right) \,.
\eeq
Using either approach to evaluate them in general and for the torus, one finds
\begin{eqnarray}
 \Gamma^r{}_{\theta\theta} 
 &=& -\frac12 (R(r)^2{})'
  = - R(r) R'(r) = (a+b\cos\chi) \sin\chi \,,\\
 \Gamma^\theta{}_{r\theta}
  &=& \Gamma^\theta{}_{\theta r} 
  = \frac{R'(r)}{R(r)} =-\frac{\sin\chi}{a+b\cos\chi}\,,
\end{eqnarray}
where we use $\chi$ and $r/b$ interchangeably as convenient.

Suppose $(x^i(\lambda))=(r(\lambda),\theta(\lambda))$ is a parametrized curve in the surface,
leading to the parametrized position vector 
\beq
  \vec r (r(\lambda),\theta(\lambda)) 
  = R(r(\lambda)) \langle \cos\theta(\lambda), \sin\theta(\lambda)  \rangle 
  +  Z(r(\lambda)) \langle 0,0,1\rangle
\,.
\eeq
The tangent vector to this curve can be evaluated using the chain rule
\beq
  \frac{d\ }{d\lambda}   \vec r (r(\lambda),\theta(\lambda))
  = \pd{\vec r(r,\theta)}{r} \frac{dr(\lambda)}{d\lambda}
   +\pd{\vec r(r,\theta)}{\theta} \frac{d\theta(\lambda)}{d\lambda} 
= \frac{dr(\lambda)}{d\lambda} e_r + \frac{d\theta(\lambda)}{d\lambda}  e_\theta\,,
\eeq
or in an abbreviated notation
\beq
  \frac{d\vec r }{d\lambda}  
  = \frac{dr}{d\lambda} e_r + \frac{d\theta}{d\lambda}  e_\theta
  =  \frac{dx^i}{d\lambda} e_i\,.
\eeq
In order to compute the derivative of this tangent vector along the curve within the surface, we need the product and chain rules. The derivative of the basis vectors along the curve, ignoring the contribution along the normal, i.e., using the covariant derivative, is the chain rule application
\beq
  \frac{d e_i}{d \lambda} 
  = \frac{\partial e_i}{\partial x_j} \frac{dx^j}{d\lambda} 
\rightarrow
  \frac{D e_i}{d \lambda} 
  = \Gamma^k{}_{ji} \frac{dx^j}{d\lambda} e_k
\eeq
and defines the covariant derivative of the basis tangent vectors along the curve itself.
The covariant derivative of the tangent vector along the curve is then a product/sum rule application to the following linear combination of products
\begin{eqnarray}
    \frac{D^2x^i}{d\lambda^2} e_i
    &\equiv& \frac{D}{d \lambda} \left(\od{x^i}{\lambda} e_i\right)
    = \frac{d^2 x^i}{d\lambda^2} e_i + \frac{dx^i}{d\lambda} \frac{D e_i}{d \lambda} 
    = \frac{d^2 x^i}{d\lambda^2} e_i + \frac{dx^i}{d\lambda} \Gamma^k{}_{ji} \frac{dx^j}{d\lambda} e_k
\nonumber\\
 &=& \left( \frac{d^2 x^i}{d\lambda^2} + \Gamma^i{}_{jk} \od{x^j}{\lambda} \od{x^k}{\lambda}
   \right) e_i\,,
\end{eqnarray}
after a convenient relabeling of the indices which are summed over.

A geodesic parametrized by an affine parameter by definition satisfies the following system of ordinary differential equations
\beq
 \ood{x^i}{\lambda} + \Gamma^i{}_{jk} \od{x^j}{\lambda}\od{x^k}{\lambda} = 0\,,\qquad
\mbox{i.e.,}\qquad
 \ood{x^i}{\lambda} + \sum^2_{j=1}   \sum^2_{k=1} \Gamma^i{}_{jk} \od{x^j}{\lambda}\od{x^k}{\lambda} = 0 \,,
\eeq
for $i=1,2$ recalling that we use $i,j,k=1,2$ interchangeably with $i,j,k = r,\theta$.
These equations are simply a statement that the tangent vector to the curve does not change at all within the surface as one moves along the curve. 

One can easily see that any linear change of parametrization $\lambda=a\bar\lambda+b$ where $a$ and $b$ are constants will preserve this form of the differential equations. It turns out that this one-parameter family of parametrizations of a given geodesic is the only freedom left in the choice of parameter in order that the tangent vector have constant length
\beq
  \frac{d\vec r}{d\lambda} \cdot \frac{d\vec r}{d\lambda}
      = \frac{ds^2}{d\lambda^2} = \left(  \frac{ds}{d\lambda} \right)^2
      = g_{ij}  \od{x^i}{\lambda} \od{x^j}{\lambda} = \mbox{\it const}\,.
\eeq
In particular one can always introduce an arc length parametrization along the geodesic by choosing $ds/d\lambda=1$.
Thus an affine parametrization is simply any parametrization for which the parameter is a linear function of some arc length parametrization
$ \lambda = A s+ B$, where $A,B$ are constants.

In the present case these geodesic equations are explicitly
\begin{eqnarray}
\ood{r}{\lambda} +\Gamma^r{}_{\theta\theta} \left(\od{\theta}{\lambda}\right)^2
&=& \ood{r}{\lambda} - R'(r) R(r) \left(\od{\theta}{\lambda}\right)^2 
  =0 \,,
\\
\ood{\theta}{\lambda} +2\Gamma^\theta{}_{r\theta} \od{r}{\lambda}\od{\theta}{\lambda}
&=& \ood{\theta}{\lambda} +2\frac{R'(r)}{R(r)} \od{r}{\lambda}\od{\theta}{\lambda}
\nonumber\\
 &=& R(r)^{-2} \odd{\lambda}\left( R(r)^2 \od{\theta}{\lambda}\right)
  =0 \,,
\end{eqnarray}
where the final form of the second equation is easily verified by expanding the derivative.
One needs to specify the initial position and initial tangent vector in order to determine a unique geodesic through a given point on the surface, aimed in a particular direction tangent to the surface. This is the geodesic initial value problem, appropriate for this coupled system of two second order nonlinear differential equations. If one chooses $\lambda=0$ for the initial parameter value, the length of the initial tangent vector determines the scaling of $\lambda$ with respect to arc length along the resulting geodesic.
It is also natural to consider geodesics between two fixed points on the surface, the two point geodesic boundary value problem, but this usually has nonunique solutions, if not an infinite number, and numerical solution requires much more sophisticated approximation schemes compared to the initial value problem. Finding a shortest such geodesic is also a nontrivial matter.

We can study these geodesic equations analytically (when possible) and numerically (when not). 
However, one can also introduce some additional quantities which enable one to more easily interpret how solutions of this system of differential equations behave. For those who are already familiar with differential geometry, this can be the starting point for the remaining discussion.

Before moving on we can note some simple solutions to these geodesic equations without specifying $R(r)$. For example, the $\theta=\theta_0$ coordinate circles (the meridians on the torus) have $d\theta/d\lambda=0=d^2\theta/d\lambda^2$, which satisfies the angular differential equation and reduces the radial equation to $d^2r/d\lambda^2=0$, with solution $r=c_1\lambda+c_2$, so that $r$ is an affine parameter on these special geodesics (recall it is actually an arc length parameter on its coordinate lines). Similarly for those points $r=r_0$ for which $R'(r_0)=0$, the radial equation is automatically satisfied while $\theta$ becomes an affine parameter on these circles of revolution. For the torus, $0=R'(r)\propto \sin(\chi)$ leads to the values $\chi=0,\pm \pi$ describing the outer and inner equators, which are therefore geodesics. For the ordinary sphere centered at the origin, the special geodesics of these two types correspond to the lines of longitude and the equator. For a horizontal plane seen as a horizontal line through the $z$-axis revolved around that axis, we only get the straight lines through its origin (which is the intersection of the plane with the $z$-axis) as these special geodesics, namely meridians.

\section{The physics approach}

So far this is the typical mathematical approach. However, let's adopt the point of view of a physicist, imagining tracing out a geodesic by identifying the affine parameter $\lambda$ with the time, so that the mental picture is now of a point particle moving on the surface, tracing out a path called the orbit of the particle, which refers to the unparametrized curve.
The tangent vector of the geodesic curve is now called the velocity
\beq
  \vec v = \frac{d\vec r}{d\lambda}
         = v^r e_r + v^\theta e_\theta
\eeq
and we might as well adopt a 2-component vector notation for components with respect to the basis vectors $e_i$ 
\beq
v^i = \od{x^i}{\lambda}:\qquad
\langle v^r,v^\theta\rangle 
= \left\langle \od{r}{\lambda},\od{\theta}{\lambda} \right\rangle
\eeq
and its magnitude $v = (\vec v \cdot \vec v)^{1/2}$ is the speed, which is just the time rate of change of the arc length along the curve (time rate of change of distance traveled)
\beq
   v = \left( g_{ij} \od{x^i}{\lambda}\od{x^j}{\lambda}  \right)^{1/2} 
     = \od{s}{\lambda} \,.
\eeq
Note that $v^{\hat r}=v^r$ is just the radial velocity while $v^\theta$ is the azimuthal angular velocity and $v^{\hat\theta}= R(r)v^\theta$ is the azimuthal component of the velocity vector. The velocity can be represented in terms of polar coordinates in the tangent plane to make explicit its magnitude and inclination angle with respect to the radial direction within the surface. To introduce these, we simply represent the orthonormal components in terms of the usual polar coordinate variables in this velocity plane in which $v^{\hat r}$ is along the first axis and $v^{\hat\theta}$ is along the second axis
\beq
  \langle v^{\hat r},v^{\hat\theta} \rangle 
  = v  \langle \cos\beta,\sin\beta \rangle \,.
\eeq
The speed plays the role of the radial variable in this velocity plane, while the angle $\beta$ gives the direction of the velocity with respect to the direction $e_{\hat r}$ in the counterclockwise sense in this plane. For affinely parametrized geodesics, the speed is constant along the geodesic; for arclength parametrized geodesics, the speed is 1.

The second order differential equations for the geodesics
\beq
  \frac{dv^{ r}}{d\lambda} = \cdots\,,\
  \frac{dv^{\theta}}{d\lambda} = \cdots
\eeq
may be interpreted as defining the radial and azimuthal angular accelerations necessary to keep the moving point on a path which doesn't change direction within the surface. As a system of second order differential equations, appropriate initial data consists of the initial values of the unknowns and their first derivatives, i.e., the initial position and initial velocity. For a ring torus, every geodesic except the inner equator must pass through the outer equator (as we will see below), and because of the rotational symmetry about the vertical axis, we might as well fix the initial position to lie on the outer equator $r=0$ at $\theta=0$. It then remains to specify the initial velocity at that point, but if we assume an arc length parametrization, then the velocity vector must be a unit vector, so only its direction within the vertical tangent plane at that initial position remains to be specified. Thus we have a 1-parameter family of arc length parametrized geodesics which interpolate between the outer equator  (horizontal initial direction) and the meridian (vertical initial direction). A computer algebra system can be easily programmed to provide a procedure for numerically solving the differential equations with these initial conditions as a function of the initial direction specified by the initial value $\beta_0$ of the angle $\beta$, which at this initial position is the direction measured from the vertical direction in the clockwise direction as seen from outside the torus looking towards the surface.

\begin{figure}[t] 
\typeout{*** EPS figure torus2: initial data}
\begin{center}
\includegraphics[scale=0.3]{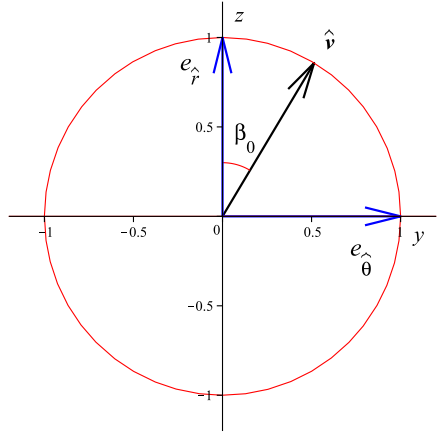}
\end{center}
\caption{
The tangent plane $x=3$ at the initial position $\langle 3,0,0\rangle$ at the intersection of the prime meridian and the outer equator of the unit ring torus where the tangent plane is vertical. The initial direction of a geodesic (unit velocity vector $\hat v=\langle \cos\beta_0,\sin\beta_0\rangle$) is here parametrized by the angle $\beta_0$ measured from the upward vertical axis in the clockwise direction as seen from the outside looking in.
} 
\label{fig:torus2}
\end{figure}

Using such a computer algebra system program one can visualize the geodesics which result from various angles $\beta_0$. Decreasing the angle from $\beta_0=\pi/2$ towards zero, one sees that the path rises higher and higher  before falling back towards the outer equator until at a certain point it goes over the North Polar Circle. Decreasing the angle further, the path creeps closer and closer to the inner equator, and then at a certain point, the path breaks through and ``threads the donut hole," i.e., crosses the inner equator.
This behavior is easily quantified as we will see below. 

However, an interesting class of geodesics are the closed ones which return to the initial position after an integral number $n\ge0$ of revolutions around the $z$-axis while making $m\ge0$ vertical oscillations, either having a point in common with the inner equator ($p=1$) or not ($p=0$), and thereafter retrace the same path over and over. 
While the Cartesian coordinates $(x,y,z)$ are periodic functions of $\lambda$ for all geodesics, these three functions only have a common period for the closed such geodesics, i.e., the position vector $\vec r$ is a periodic vector function of $\lambda$.
We will see below that if a geodesic returns to the initial point, its velocity must have either the same initial angle or its complement, and in the latter case, it will return to the initial point with the same initial angle after twice as many azimuthal revolutions, which follows from reflection symmetry across the outer equator.
These closed geodesics can therefore be classified by their period number pairs and the inner equator parameter: $[m,n;0]$ or $[m,n;1]$. By trial and error one can find some of these special geodesics using a numerical program, making it a sort of mathematical computer game. For some reason this fascinated me and perhaps it might do the same for some students. A Maple worksheet is available to test this out at my website. However, why guess if one can predict? This is possible as described below.

\section{Energy and angular momentum: the key to unlocking the mystery}

To understand the system of two second order geodesic equations, one can use a standard physics technique of partially integrating them and so reduce them to two first order equations by using two constants of the motion which arise from the two independent symmetries of the ``equations of motion": time translation, which results in a constant energy, and rotational symmetry, which results in a constant angular momentum. There is no need to understand how this happens---it is a topic which belongs in an advanced mechanics course. Here we will just manipulate our specific equations to accomplish our goals.

Since the mass $m$ of the point particle whose motion traces out a geodesic path is irrelevant in this problem, those physical quantities like energy and momentum which involve the mass as a proportionality factor will instead by replaced by the ```specific'' quantities obtained by dividing out the mass.
Thus the specific kinetic energy (kinetic energy $\frac12 m v^2$ divided by the mass $m$) is just
\begin{eqnarray}
  E &=& \frac12 v^2 = \frac12 g_{ij} \od{x^i}{\lambda}\od{x^j}{\lambda}
\\
   &=& \frac12 \left(\od{r}{\lambda}\right)^2 
    + \frac12 R(r)^2\left(\od{\theta}{\lambda}\right)^2
\\
   &= &\frac12 v^2 \cos^2\beta
     +\frac12 v^2 \sin^2\beta \,,
\end{eqnarray}
and vice versa
\beq 
  v = (2E)^{1/2} \,.
\eeq
Both the specific energy and speed are constant for an affine parametrization of the geodesic.

In the physics approach the specific energy of the particle is constant because from the point of view of its motion in space, it is only accelerated perpendicular to the surface (since the component of the second derivative along the surface has been designed to be zero). If a force were responsible for this acceleration, namely the normal force which keeps the particle on the surface, since it is perpendicular to the velocity of the particle it would not do any work on the particle (recall that work done is the integral of the dot product of the force and velocity vectors). Thus its energy  and hence specific energy $E$ must remain constant (said to be ``conserved'' along the geodesic). Similarly the speed $ v = (2E)^{1/2}$ must therefore be constant along a geodesic according to this reasoning.

In this same physics language, the second geodesic equation tells us that the specific angular momentum about the axis of symmetry (defined exactly as in the case of circular motion around an axis with radius $R(r)$, namely the velocity $v^{\hat\theta}= R(r)\, d\theta/d\lambda$ in the angular direction multiplied by the radius $R(r)$ of the circle)
\begin{eqnarray}
  \ell &=& \hat z \cdot (\vec r \times \vec v)
\nonumber\\
&=& R(r)\, v^{\hat\theta}
       = R(r)^2 \od{\theta}{\lambda}  
       = R(r)\, v \sin\beta
\end{eqnarray}
is constant along a geodesic, or ``is a conserved quantity" for the motion.  This relation can then be used to re-express the variable angular velocity $d\theta/d\lambda$
in the specific energy  formula in terms of the constant angular momentum
to obtain the constant specific energy entirely expressed in terms of the radial motion (i.e., only involving the variables $r$ and $dr/d\lambda$) and another constant of the motion
\beq
 E = \frac12 \left(\od{r}{\lambda}\right)^2 
    + \frac12 \frac{\ell^2}{R(r)^2} \,.
\eeq
In other words we have reduced the problem to motion in one dimension governed only by a first order differential equation for the radial coordinate. Before continuing its description in the next section, some general considerations are worth mentioning. For simplicity we drop the ``specific" qualifier from the energy and angular momentum. Note also that for a ring torus where the azimuthal radius $R(r)$ cannot vanish, if the angular momentum $\ell$  is zero, then $d\theta/d\lambda$ must vanish leading to a constant value of the azimuthal angle $\theta=\theta_0$ corresponding to purely radial motion. This is a trivial geodesic as shown above. We can therefore concentrate our attention on the nonradial geodesics.

\begin{figure}[t] 
\typeout{*** EPS figure torus initial velocity plane}
\begin{center}
\includegraphics[scale=0.35]{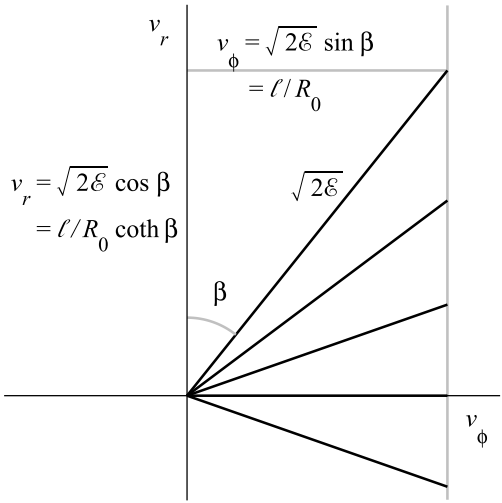}
\end{center}
\caption{
Instead of fixing the speed/energy (circles) in the initial velocity plane at the outer equator, fixing the angular momentum leads to a vertical line corresponding to all geodesics in either counterclockwise ($\ell>0$) or clockwise ($\ell<0$) azimuthal motion or in purely radial motion ($\ell=0$). Here the counterclockwise case is illustrated.
} 
\label{fig:torusvelocityplane}
\end{figure}

Re-expressing the velocity in terms of the two constants of the motion yields
\beq
    \langle v^{\hat r},v^{\hat\theta} \rangle 
  =   \langle (2E)^{1/2}\cos\beta,(2E)^{1/2}\sin\beta \rangle
  =    \left \langle \od{r}{\lambda},\frac{\ell}{R(r)} \right\rangle \,.
\eeq
The second component of this equation gives a direct relationship between the angle of inclination and the  energy for a given  angular momentum for nonradial motion $\beta\ne0$, $\ell\ne0$
\beq\label{eq:ellE}
 \frac{\ell}{(2E)^{1/2}} = R(r) \sin\beta
\eeq
or equivalently
\beq
  E = \frac{\ell^2}{2 R(r)^2\sin^2\beta}
    = \frac{\ell^2}{2 R(0)^2\sin^2\beta_0}
 \,.
\eeq
The latter relation is useful to  geometrically interpret the initial velocity data at $r=0$ for this problem.
When $E=\frac12=\frac12 v^2$, the speed $v$ is 1 which corresponds to an arc length parametrization.
However, we will see soon that it is more convenient to fix the angular momentum than the energy. For example if we choose instead $\ell=1$, then varying the initial angle $\beta_0$ varies the  energy, as illustrated in Fig.~\ref{fig:torusvelocityplane}. Different energy levels then describe the choice of initial direction relative to the $r$ coordinate circles (meridians).

To interpret the former relation, note that
dividing the specific angular momentum (units of length times speed) by the speed leads to the specific angular momentum in a unit speed (arc length) parametrization, so that its units become just those of length:  
$\ell/(2E)^{1/2} = R(r) \sin\beta$. In fact this is just the vertical component of the unit speed angular momentum vector
\beq\label{eq:unitspeedangmom}
    \ell/(2E)^{1/2} = \hat z \cdot (\vec r \times \hat v) =R(r)\sin\beta \,.
\eeq
This combination of the constants of the motion is of course also constant along a geodesic,
so that for initial data at $(r,\theta)=(0,0)$ one has
\beq
R(r)\sin\beta = R(0)\sin\beta_0 \,.
\eeq
The obvious interpretation of this constant length  is its value at a radial ``turning point" $r=r_{\rm(ext)}$ where $\sin\beta=\pm1$ (see below), if it exists, where the unit velocity is aligned with the azimuthal angular direction and this reduces to $\pm R(r_{\rm(ext)})$.
The existence of this constant is a consequence of the 1-parameter rotational group of symmetries of the torus (or indeed any surface of revolution) about its axis of symmetry, and indeed such a constant of the motion arises when the surface is invariant under any 1-parameter group of symmetries, a fact easily seen in the variational approach to the geodesic equations. In the mathematical language, this quantity is a constant discovered by Clairaut for geodesic motion in a surface described in a coordinate system adapted to this 1-parameter group of symmetries \cite{oprea,pressley}. However, it appears that none of the textbooks on differential geometry develop the tools of transformation groups or symmetry groups (associated with the names of Lie and Killing), and so cannot remark on the interpretation available in most elementary textbooks on general relativity. These surfaces of revolution are invariant under the the 1-parameter group of rotations about the symmetry axis, and as such are ``isometries" of the intrinsic metric. The circular orbits of this symmetry group are the integral curves of a Killing vector field generator, which in this case is $\partial/\partial \theta$. Translations in this coordinate $\theta$ correspond to these rotations. The Clairaut conditions on Clairaut coordinates just state that the coordinates are orthogonal and one is a Killing coordinate, i.e., such that the metric is invariant under translation in that coordinate, and hence the components of the metric in the coordinate system adapted to that single independent Killing vector field are independent of that coordinate. The constant of the motion is just the inner product of the unit tangent vector (unit velocity) along a geodesic with that Killing vector field
\beq
  (\hat v)_\theta= \hat v \cdot \frac{\partial}{\partial\theta} = R(r) \sin(\beta) \,.
\eeq

We can use this constant of the motion to answer some interesting questions about the crossing angle of a geodesic with the parallels it encounters. For example,
if a geodesic starting at the outer equator crosses the inner equator, it must do so at an angle which satisfies
\beq
 \sin\beta(\lambda) = \frac{R(0)}{R(b\pi)} \sin\beta_0 
  =\frac{a+b}{a-b} \sin\beta_0 
  =\frac{2+c}{c} \sin\beta_0 
\,.
\eeq
This crossing angle $\beta(\lambda)$ is a bigger angle than $\beta_0$ for a ring torus for which $c>0$.
If a geodesic starts out at the outer equator and returns to the outer equator $R(r(\lambda))=R(0)$, then $\sin(\beta(\lambda))=\sin(\beta_0)$, so its crossing angle with respect to the outer equator is either $a)$ the same $\beta_0$ or $b)$ its supplement $\pi-\beta_0$. This crossing angle cannot equal one of the two sign-reversed angles $-\beta_0,-(\pi-\beta_0)$ since the azimuthal angle $\theta$ is a monotonically increasing or decreasing function along a geodesic for nonzero angular momentum (the sign of $d\theta/d\lambda$ does not change), always revolving either clockwise or counterclockwise around its symmetry axis as seen from above. If the geodesic returns to the initial starting point at the origin of coordinates, the first case $a)$ means the initial data at the final point is the same as at the starting point and so it must retrace its path and the geodesic is closed, while in the second case $b)$ by symmetry with respect to reflections across the equatorial plane, the curve must be a geodesic which will return to the starting point with the original initial velocity after another time interval of equal length $\lambda$, and is therefore also closed. However, for geodesics which pass through the inner equator, the radial coordinate is either monotonically increasing or decreasing ($dr/d\lambda$ always has the same sign, never vanishing), and so cannot cross itself at the starting point with the supplementary angle and hence must begin retracing its path as a closed geodesic. In fact since for such geodesics both $r$ and $\theta$ are monotonic functions of $\lambda$, the velocity must always remain within one of the four quadrants of the tangent space for which the sign combinations of $(dr/d\lambda,d\theta/d\lambda)$ are fixed. Thus the unbound geodesics either never intersect themselves or they are closed and retrace their own path. The bound geodesics may intersect themselves, but only at supplementary angles with respect to the meridians.

The so called ``turning points" of the radial motion occur at extreme radial values $r=r_{\rm(ext)}$, where the radial velocity $v^{\hat r}=v\cos\beta$ is zero and the velocity purely azimuthal. These points correspond to $\cos\beta=0$ and so $\sin\beta=\pm1$ and are therefore determined as the roots of the equation
\beq
  R(r_{\rm(ext)}) = R(0)|\sin\beta_0| \,,
\eeq
or explicitly for the torus,   $ r_{\rm(ext)} = \pm r_{\rm(max)}$ where
\begin{eqnarray}\label{eq:rmax}
   \frac{r_{\rm (max)}}{b}   =\chi_{\rm(max)}
&=& \arccos\left( \frac{(a+b)|\sin\beta_0|-a}{b} \right)
\nonumber\\
    &=& \arccos\left( (2+c)|\sin\beta_0|-(1+c) \right) 
\,.
\end{eqnarray}
This can only hold for the so called bound geodesics (see below) for which a turning point radius $R(r_{\rm(ext)})$ exists and hence must be larger than or equal to the inner equatorial radius
\beq
  R(0)|\sin\beta_0|  =   R(r_{\rm(ext)}) \ge R(b\pi) \leftrightarrow
   -\pi b \le r_{\rm(ext)} \le \pi b
\,,
\eeq
so that $\sin\beta_0$ is larger in absolute value than the ratio of the inner and outer equatorial radii
\beq
|\sin\beta_0|  \ge    \frac{R(b\pi)}{R(0)}
         = \frac{a-b}{a+b} =\frac{c}{2+c}
\,.
\eeq

The formula for $r_{\rm(ext)}$ can be used to aim geodesics starting at the outer equator to hit a certain extremal radius. 
For example, to asymptotically approach the inner equator one needs the initial angle to satisfy
\beq
|\sin\beta_0| = \frac{a-b}{a+b} = \frac{c}{2+c} \rightarrow \beta_0 = \pm\beta_{\rm(crit)}, \pi \pm\beta_{\rm(crit)}
\eeq
where
\beq
   \beta_{\rm(crit)} = \arcsin \left(\frac{c}{2+c}\right)
\eeq
($19.5^\circ$ for the unit ring torus) defines a critical angle,
or to reach the Northern or Southern Polar Circles the angle must be larger and satisfy
\beq
|\sin\beta_0| = \frac{a}{a+b}=\frac{1+c}{2+c}
\eeq
defining another special angle
\beq
   \beta_{\rm(polar)} = \arcsin \left(\frac{1+c}{2+c}\right)
\eeq
 ($41.8^\circ$ for the unit ring torus). All initial angles smaller than this correspond to geodesics which cross the polar circles, while all initial angles smaller than the former critical angle result in unbound geodesics which cross the inner equator.

Note that one can simply solve the energy equation describing the radial motion for the radial velocity and integrate
\beq\label{eq:lambdaofr}
  \lambda-\lambda_0 = \int_{r_0}^r \frac{1}{(2E-\ell^2/R(r)^2)^{1/2}} \, dr\,.
\eeq
For the torus this integral can be done exactly in terms of elliptic functions, but the result cannot be inverted to express $r$ as a function of $\lambda$. However, it can also just be integrated numerically to provide $\lambda(r)$, which becomes arc length once multiplied by the speed. Note also that at turning points of the radial motion $dr/d\lambda=0$, so that the denominator of this integral vanishes. Thus an integral for which the turning point $r$ is an upper limit of integration is an improper integral, which even if it converges can still lead to some complications for numerical evaluation.

\section{The reduced problem of the radial motion}

\begin{figure}[t] 
\typeout{*** EPS figure torus potential}
\begin{center}
\includegraphics[scale=0.4]{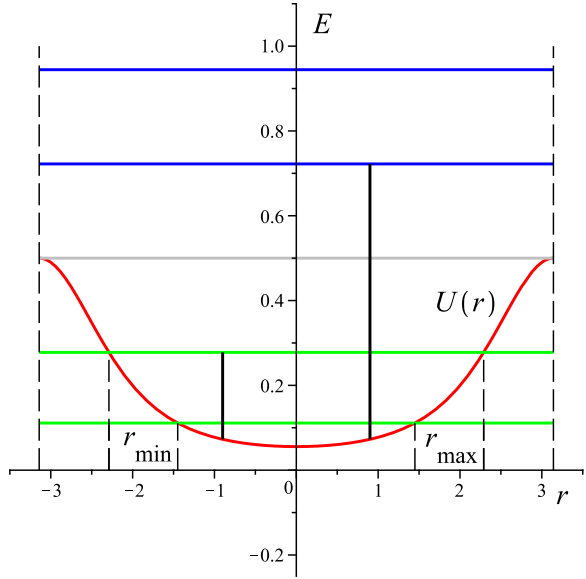}
\end{center}
\caption{
The unit  angular momentum potential (i.e., for $\ell^2=1$) for ``nonradial" torus geodesics starting on the outer equator on the unit ring torus $(a,b)=(2,1)$, where $r=\chi$. If the energy is less than 1/2, the geodesics do not make it through the hole: $|r_{\rm max}|<\pi$. Of course the potential is periodic as well, and we have only shown one period here: $-\pi\le r\le \pi$. For energy larger than 1/2, the variable $r$ is unbounded. Five energy levels are shown: 2 typical bound geodesics ($E<1/2$), 1  unique bound geodesic asymptotically approaching the inner equator ($E=1/2$), and 2 typical unbound geodesics ($E>1/2$). The bound geodesics have symmetric extremal radii. The potential minimum corresponds to the outer equator, and the potential maximum to the inner equator.
} 
\label{fig:toruspotential}
\end{figure}

\begin{figure}[t] 
\typeout{*** EPS figure torus potential unbound periodic [m,1;1]}
\begin{center}
\includegraphics[scale=0.4]{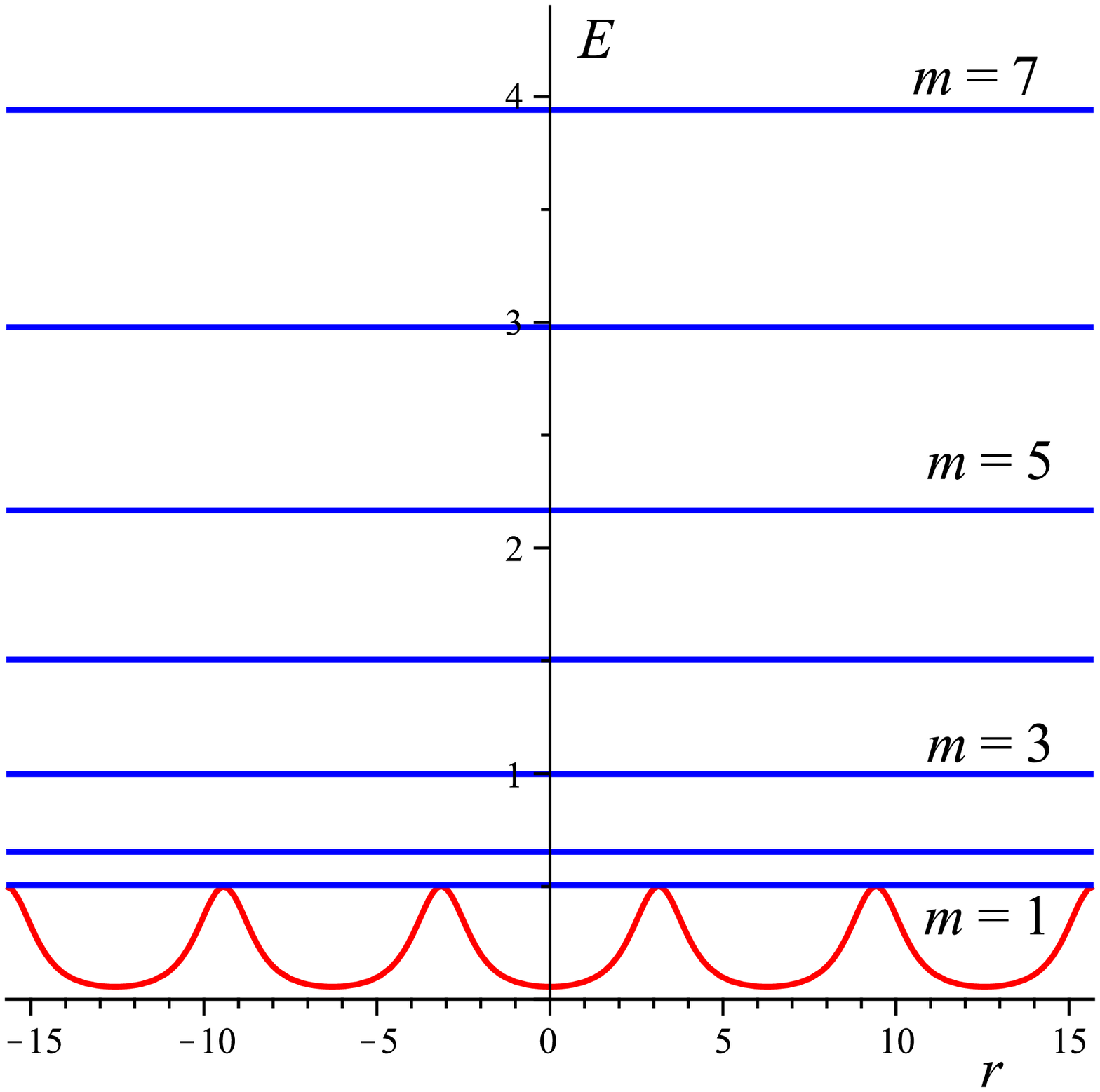}
\end{center}
\caption{
The energy levels $E>\frac12$ for the unbound geodesics on the unit ring torus of type $[m,1;1]$ for $m=1,\ldots,7$. Here the periodic potential merely slows down or speeds up their monotonic unbounded radial motion.
} 
\label{fig:toruspotentialperiodic}
\end{figure}

\begin{figure}[t] 
\typeout{*** EPS figure torus periodic geos}
\vglue-1in
\begin{center}
\includegraphics[scale=0.7]{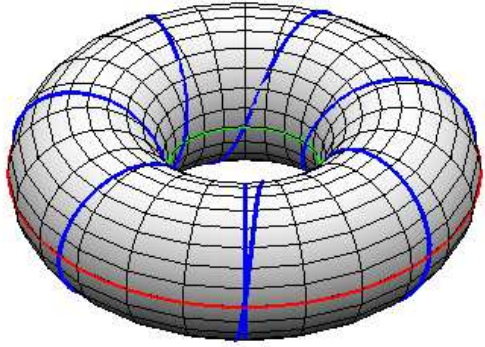}
\end{center}
\vglue-1in
\caption{
An closed geodesic $[7,1;1]$ for the unit ring torus with 7 loops crossing the inner equator, corresponding to an initial angle of $\beta\approx 0.119$ ($\approx 6.81^\circ$) and relatively high energy $E=3.942$ on the unit ring torus.} 
\label{fig:torusperiodicgeos}
\end{figure}

The real advantage of our physics approach is in the diagram we can create to describe the radial motion. The second term in the energy is the ``effective potential energy'' due to the angular motion, called the centrifugal potential. If we only consider nonradial motion, we can fix a convenient value of the angular momentum while allowing the energy to be determined by the initial angle and therefore have a single potential energy function $U$ describe the radial motion
\beq
   \frac12 \left(\od{r}{\lambda}\right)^2 
    + U = E \,,\qquad
  U = \frac12 \frac{\ell^2}{R(r)^2} = \frac12 \frac{\ell^2}{(a+b\cos(r/b))^2}  \,.
\eeq
Since the sum of the radial kinetic energy and this potential energy function is the constant energy for a given geodesic, once one sets an energy level, the motion must be restricted to the region where the difference $E-U$ is nonnegative, with the radial velocity $dr/d\lambda$ necessarily zero at points where $E=U$, which determine the extreme values $r_{\rm(ext)}$ of the radial variable within which the geodesic is confined. Fig.~\ref{fig:toruspotential} illustrates one period of the periodic potential for the unit ring torus $(a,b)=(2,1)$ with the choice $\ell=1$, where the inner equator has energy $1/2$ and is the unstable equilibrium solution sitting at the peak of the potential energy graph. Geodesics with energy less than $1/2$ are bound with radial motion confined between the symmetric extreme values of the radial coordinate, analogous to the negative energy elliptic orbits of the Kepler problem. The geodesics with energy $E>1/2$ are unbound, requiring an infinite interval $-\infty < r<\infty$ of radial coordinate values to describe continuously, and therefore the periodic potential shown in Fig.~\ref{fig:toruspotentialperiodic}. For example, the energy $E=3.942$ corresponds to
the 7 loop unbound geodesic shown in Fig.~\ref{fig:torusperiodicgeos} for the initial angle $\beta_0= 0.119$ radians ($6.8^\circ$).

If we return to the original second order equation for the radial variable and replace the angular velocity using conservation of angular momentum we get
\beq
  \ood{r}{\lambda}
 = R'(r) R(r) \left(\od{\theta}{\lambda}\right)^2 
 = R'(r) R(r) \left( \frac{\ell}{R(r)^2} \right)^2
 = \frac{\ell^2 R'(r)}{R(r)^3}
 = -\frac{dU}{dr}
 \,.
\eeq
This shows that the effective potential really does generate the radial component of the effective (specific) force through its derivative. However, we are spared having to deal with this second order equation since we have its first integral already, namely the energy relation. Note that this effective force is really just the radial coordinate acceleration required to keep the geodesic direction unchanged and is not a true force. The only real force on the motion (if this were a real point particle in motion) is the force normal to the surface required to keep the curve from leaving the surface, a force which arises from the constraint that the motion remain within the surface.
For completeness we show that the constancy of energy implies the second order radial equation by a simple derivative calculation. If the energy is constant but $r$ is not constant, then
\beq
 E = \frac12 \left(\od{r}{\lambda}\right)^2 
    + \frac12 \frac{\ell^2}{R(r)^2} \,, \quad
  \frac{dE}{dt} 
   = \od{r}{\lambda} \left(\ood{r}{\lambda} 
     -\frac{\ell^2 R'(r)}{R(r)^3} \right) 
 \,.
\eeq
the expression in parentheses on the right hand side of the last equation must be zero, which is the second order radial equation.

When the potential has a critical point where the derivative is zero, we obtain an equilibrium solution for which the radial variable is constant, namely a $\theta$ coordinate line (a parallel) which is a horizontal geodesic circle of revolution. If the second derivative is positive (a local minimum), it is a stable equilibrium since for slightly larger energies the radial motion will be confined to a narrow interval of radii around it, while if the second derivative is negative, it will correspond to an unstable equilibrium since the radial motion will go far from the initial value. The torus  has a periodic potential (we must consider $r$ unbounded to describe motion which wraps multiple times around the circular cross-section) in the shape of a repeated ``potential well'' as shown in Fig.~\ref{fig:toruspotentialperiodic} for a ring torus. The bottom of the well where 
\beq
U = U(0)
   = \frac{\ell^2}{2R(0)^2} 
   = \frac{\ell^2}{2(a+b)^2}
\eeq 
is at $r=0=\chi$ at the outer equator (a stable equilibrium), while the rim of the well is at $r=\pm \pi b$ or $\chi=\pm \pi$ at the inner equator (an unstable equilibrium) where 
\beq
U  = U(b\pi)
   = \frac{\ell^2}{R(b\pi)^2} 
   = \frac{\ell^2}{2(a-b)^2} \,.
\eeq
The \lq\lq walls'' of the potential act as a barrier to reaching the inner equator unless the energy is sufficiently high to overcome this barrier. The turning points $r=r_{\rm(max)},r_{\rm(min)}$ of the radial motion for this case occur at the minimum radius $R(r)$ allowed by conservation of angular momentum: where the energy level intersects the potential graph. The maximum angular velocity $d\theta/d\lambda$ occurs at these turning points where the radial velocity is momentarily zero, simply by conservation of angular momentum since these turning  points are closest to the symmetry axis about which the geodesic is revolving. For the horn torus, the rim of the potential well moves up to positive infinity, so that the nonradial unbound geodesics disappear.

By fixing the angular momentum, one can consider a single potential function instead of a family at the cost of allowing the variable energy to describe different geodesics. For example, if we consider an initial point at the origin of coordinates on the outer equator (see Fig.~\ref{fig:torusvelocityplane}), we only have to vary the angle of the direction of the velocity to describe the 1-parameter family of distinct geodesic paths which emanate from the initial position, and the different inclination angles correspond to the different values of the energy (equivalently of the speed). 
There are a number of obvious choices to fix $\ell\neq0$ in the potential function which only depend on the shape parameter
\beq
  U =  \frac{\ell^2}{2R(r)^2} = \frac{\ell^2}{2(a+b\cos(r/b))^2} 
    = \frac{\ell^2/b^2}{2(c+1+\cos(r/b))^2} \,.
\eeq
For example, setting $|\ell|=R(0),R(b\,\pi),b$ or equivalently
$|\ell|/b = (a+b)/b,(a-b)/b,1 = c+1,c,1$
leads to values at the outer and inner equators respectively of
\beq
  U(0)= \frac12,\frac12 \left(\frac{c}{c+2}\right)^2, \frac{1}{2(c+2)^2}\,,\
  U(b\,\pi)= \frac12 \left(\frac{c+2}{c}\right)^2,\frac12, \frac{1}{2c^2}\,.
\eeq
The middle choice won't work for horn tori where $c=0$.
For the unit ring torus ($c=1$) either of the last two choices are the same as simply setting $|\ell|=1$, which will be assumed below in the unit ring torus examples. 

This latter potential is illustrated in Figs.~\ref{fig:toruspotential} and \ref{fig:toruspotentialperiodic} for the unit ring torus.
Those orbits for which $|\chi_{\rm (max)}|\le \pi$ are referred to as the bound orbits since the range of the radial variable is confined to a finite interval of values which is symmetric about the center $r=0$ because of the reflection symmetry of the potential about its center. 
In the physics language, the bound geodesics are ``trapped" in the potential ``well," while the unbound geodesics are free to pass from one well to the next as they continue to loop around the torus forever.
For the energy level $\frac12$, there are two geodesics. One is the inner equator which is an unstable equilibrium, and the other is a geodesic which asymptotically approaches the inner equator both from above and from below at the two limiting turning points of the radial motion, revolving an infinite number of times around the symmetry axis in each direction.
If the energy is larger than $\frac12$, the orbits are referred to as unbound since one requires an unbounded range of the radial variable to describe them continuously. Note that only the inner equator geodesic does not pass through the outer equator (obviously!), so initial data at a point on the outer equator is enough to describe all other geodesics.

Thus from the potential diagram typical of all ring tori, one sees five classes of geodesics: the inner and outer equators, the geodesic asymptotic to the inner equator, the remaining bound geodesics, and the unbound geodesics, to which we must add the sixth type, namely the radial geodesics since the radial coordinate lines are geodesics for any surface of revolution, as they must be because of the reflection symmetry about the plane cross-section of the surface whose intersection is the given radial coordinate line. 

\begin{figure}[t] 
\typeout{*** EPS figure torus potential unbound periodic [m,1;1]}
\begin{center}
\includegraphics[scale=0.3]{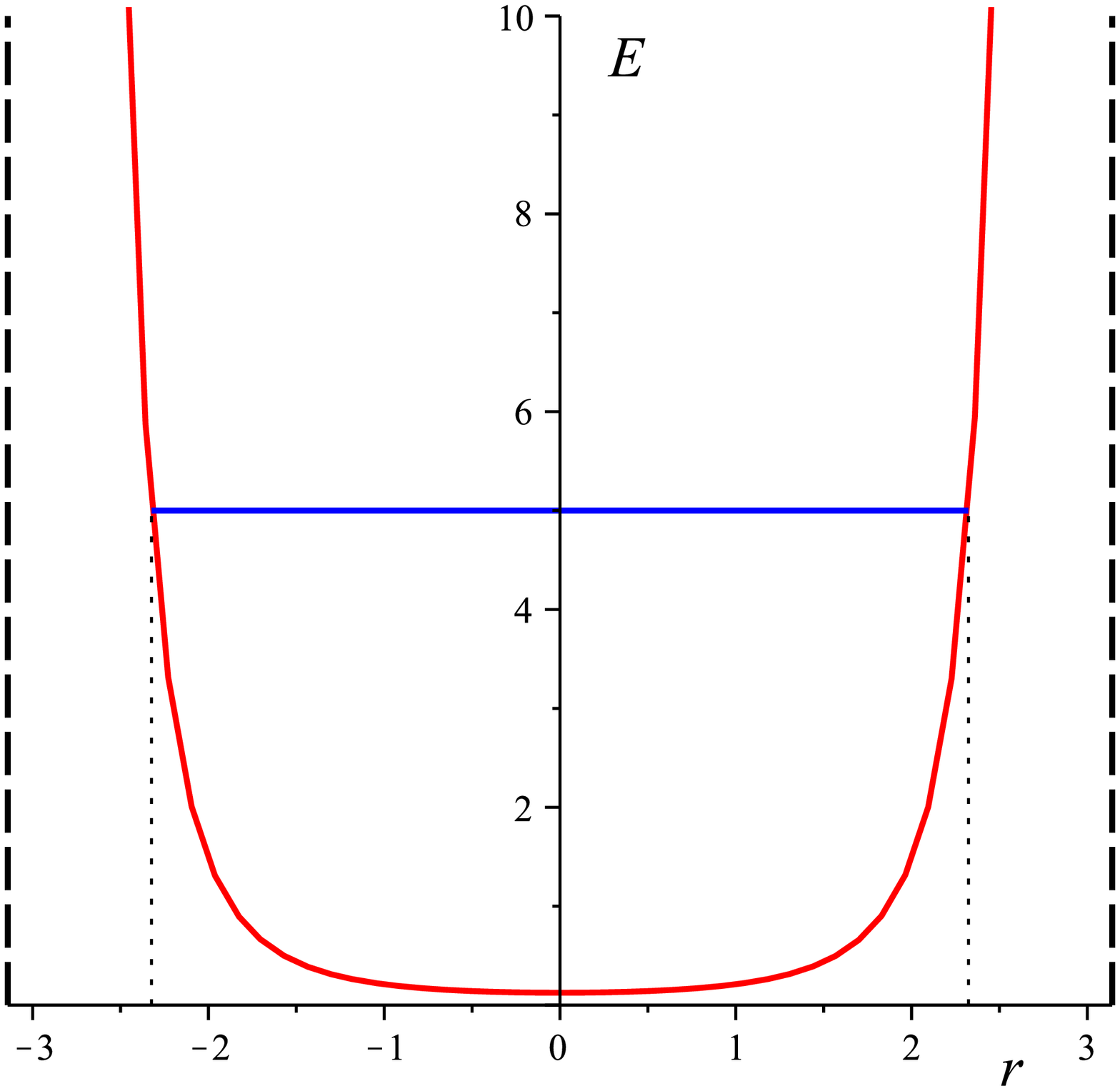}
\includegraphics[scale=0.3]{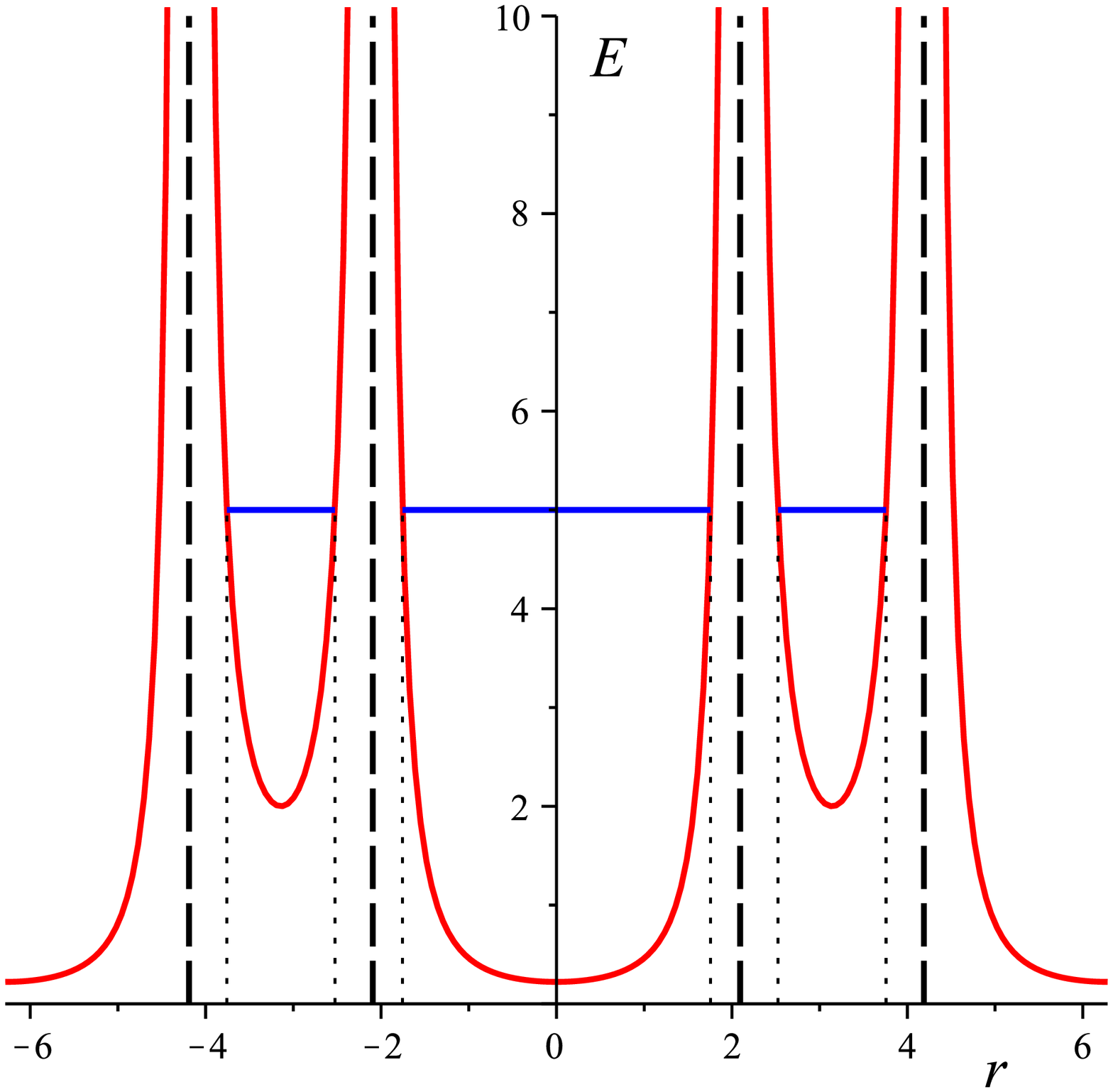}
\end{center}
\caption{
The potential for the $(a,b,c)=(1,1,0)$ unit horn torus (left) on the interval $-\pi\le r \le \pi$
and for the $(a,b,c)=(1/2,1,-1/2)$ spindle torus (right) on the interval $-2\pi\le r \le 2\pi$.
The horn torus has a single equilibrium solution at the outer equator at the origin, with infinite potential walls surrounding the inner equator.
The spindle torus has equilibrium solutions at both the outer (apple) equator (center well) and the
inner (lemon) equator (outer wells), and nonradial geodesics are confined to either the lemon or spindle. Turning points for one energy level are shown for both cases.
} 
\label{fig:toruspotentialimproper}
\end{figure}

On the other hand when one passes from ring tori to the horn torus by decreasing the shape parameter $c$ to zero, the energy level for the inner equator grows to positive infinity and the nonradial unbound geodesics disappear. All geodesics are trapped in the infinite potential well whose walls go to infinity at the inner equator where $R(r)=0$.
Continuing into the spindle tori, one has two separate sets of nonradial bound geodesics confined either to the apple or the lemon component of the torus, and the two points of self-intersection are protected by the infinite potential walls surrounding the axis of symmetry. The infinite values of the potential occur at $\pm r_\infty/b = \pm \chi_\infty$ where $R(r_\infty)=0$, namely
\beq
  \chi_\infty = \arccos\left(-a/b\right) = \arccos(-(c+1))\,.
\eeq
There are two stable equilibrium solutions, one at the inner equator on the lemon and one at the outer equator on the apple. These potentials are illustrated in Fig.~\ref{fig:toruspotentialimproper}, where the condition $|\ell|/b=1$ is assumed so that $U(0)=1/2(c+2)^{-2}<1/2$ and $U(b\pi)=1/2 c^{-2}>1/2 $. As one decreases $c$ from 0 towards $-1$ the inner well narrows and rises as $\chi_\infty\to \pi/2$ and the outer well widens and deepens, until at the limiting case $c=-1$ of the sphere, the inner and outer wells are the same shape, and only the inner well is needed to describe the nonradial geodesics, all of which are bound and closed with $-\pi/2<r<\pi/2$.
Note that if we allow $-b<a<0$ or $-2<c<-1$, then we again get the spindle tori but with $r=0$ now describing the inner equator on the lemon surface. This allows one to describe the lemon geodesics using initial data at the origin of coordinates $(r,\theta)=(0,0)$ as in the remaining cases.

Since the outer equator is a stable equilibrium, nearby geodesics will oscillate about it with a frequency that is easily calculated by approximating the bottom of the potential well by its quadratic Taylor polynomial. The equations of motion for these limiting solutions are then
\begin{eqnarray*}
r=0:&&\qquad 
    \left.\od{\theta}{\lambda}\right|_{r=0} = \frac{\ell}{(a+b)^2} 
     \rightarrow \theta= \frac{\ell}{(a+b)^2}  \lambda\,,\\
r\ll b:&&\qquad 
E=\frac12 \left(\od{r}{\lambda}\right)^2 +\frac{\ell^2}{(a+b)^2} + \frac{ \ell^2 r^2}{2b(a+b)^3}\,,
\end{eqnarray*}
having imposed the initial condition $\theta=0$ at $\lambda=0$ on the angular solution.
Defining
\beq
 \omega=\frac{|\ell|}{[b(a+b)^3]^{1/2}}\,,\quad
 \frac12 \mathcal{R}^2 = E-\frac{\ell^2}{(a+b)^2}
\eeq
the energy equation describes a simple harmonic oscillator
\beq
  \mathcal{R}^2 = \left(\od{r}{\lambda}\right)^2+ \omega^2 r^2\,.
\eeq
Integrating this for $\lambda$ as a function of $r$ (solve for $d\lambda$ and integrate), imposing the initial condition $\lambda=0$ at $r=0$ yields the result
\beq
 r = \pm \mathcal{R} \sin(\omega \lambda) = \pm \mathcal{R} \sin((c+2)^{1/2} \theta) \,.
\eeq
Thus there are $(c+2)^{1/2}$ radial oscillations during one revolution about the symmetry axis, or $\sqrt{3}\approx 1.73$ for the unit ring torus. The final equation here, having eliminated the parameter $\lambda$, directly gives the orbit equation for this case, namely the path traced out by these limiting geodesics near the outer equator. Note that the smallest positive integer value of the positive shape parameter that leads to an asymptotically closed orbit for the small oscillations at the outer equator is $c=2$ (i.e., $a/b=3$), with 2 oscillations per revolution. For a sphere (limiting value $c=-1$) this leads to one oscillation per revolution, corresponding to a great circle, which also characterizes the exact oscillations about the equator of any amplitude.
For the redundant parameter range $-2<c<-1$ for spindle tori (duplicating the range $-1<c<0$) where $r=0$ instead corresponds to the inner equator (on the lemon) of the self-intersecting torus, one finds an infinite number of values of $c$ for which the periodic limiting small amplitude oscillations have a period which is an integer multiple of the length of the inner equator.

\begin{figure}[t] 
\typeout{*** EPS figure torus periodic geo [3,2;p]}
\vglue-1cm
\begin{center}
\includegraphics[scale=0.30]{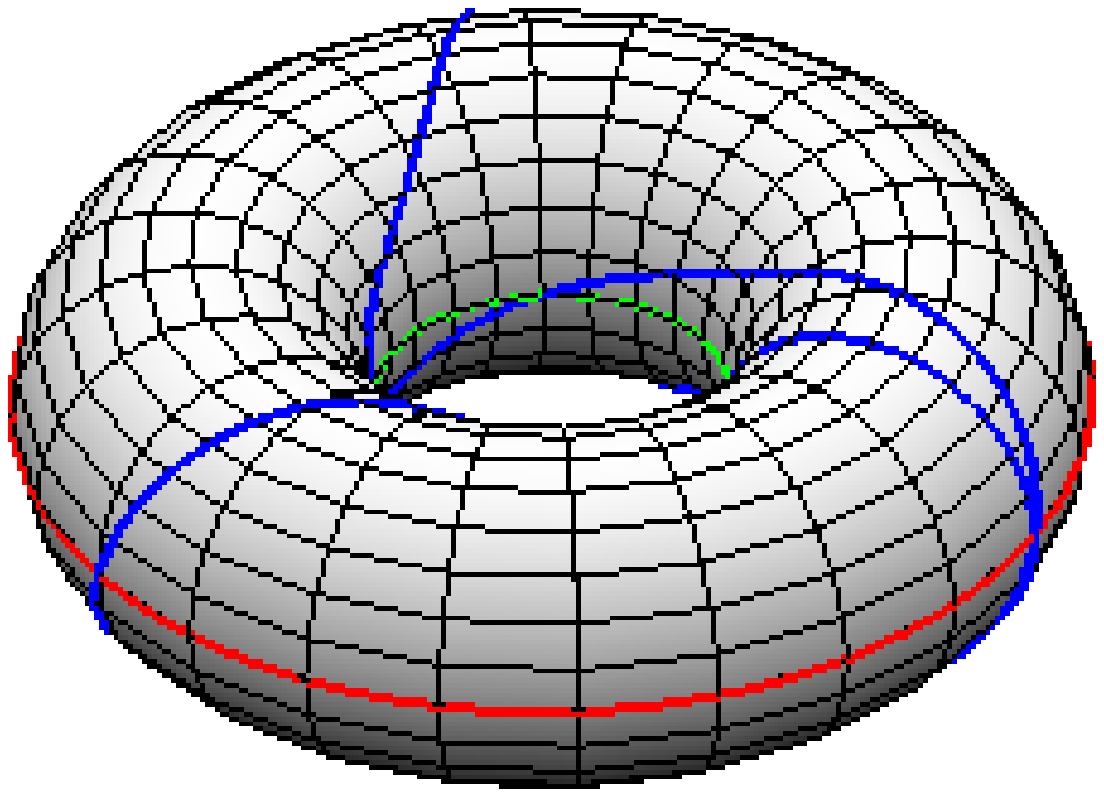}
\includegraphics[scale=0.30]{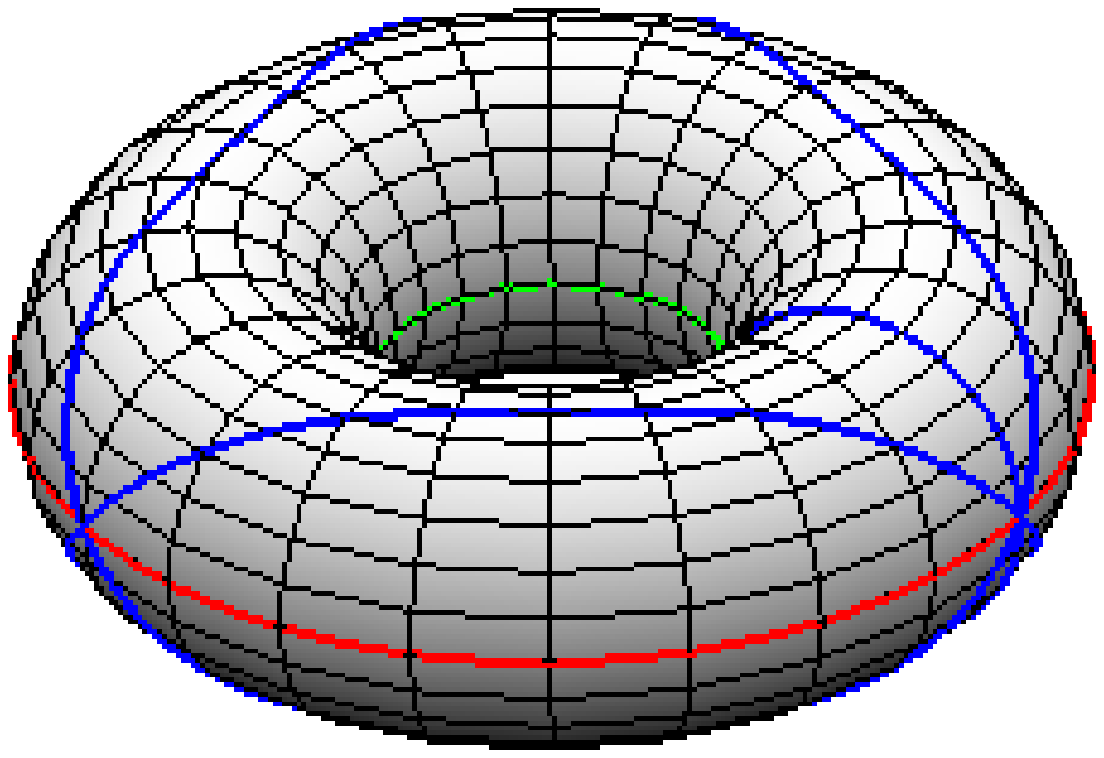}
\end{center}
\vglue-1cm
\caption{
Periodic orbits $[3,2;1]$ (left, unbound) and $[3,2;0]$ (right, bound) for the unit ring torus corresponding to an initial angles of $\beta_0=0.0.3226$ ($18.5^\circ$)  and $\beta_0=0.7167$ ($41.1^\circ$). 
}
\label{fig:torusperiodicgeos32}
\end{figure}

\section{The orbits}

If we are only interested in the paths traced out by the geodesics (the orbits), we can easily eliminate the affine parameter $\lambda$ and get a direct integral relationship between $\theta$ and $r$ in general, first solving the energy equation for the radial velocity and then using the chain rule
\beq
   \od{\theta}{r} = \od{\theta}{\lambda} / \od{r}{\lambda}
   = \frac{\ell}{R(r)^2} \frac1{(2E-\ell^2/R(r)^2)^{1/2}}
   = \frac{\ell}{R(r)} \frac1{(2E R(r)^2-\ell^2)^{1/2}}
\,.
\eeq
Note that this goes infinite at a turning point of the radial motion where $dr/d\lambda=0$.
This expression can be integrated and expressed in terms of the initial data at $(r,\theta)=(0,0)$ parametrized by the initial polar angle $\beta_0$ by using the relation
\beq
  E = \frac{\ell^2}{2 R(0)^2\sin^2\beta_0} \,,
\eeq
which leads to
\begin{eqnarray}
\theta &=& \int_0^r \frac{R(0) \sin\beta_0}{R(r)(R(r)^2-R(0)^2\sin^2\beta_0)^{1/2}} \,dr
\equiv F(r,\beta_0)
\nonumber\\
&=&  \int_0^\chi \frac{(c+1) \sin\beta_0}{(c+1+\cos\chi)[(c+1+\cos\chi)^2-(c+1)^2\sin^2\beta_0]^{1/2}} \,d\chi
\nonumber\\
&\equiv& G(\chi,\beta_0)
\,.
\end{eqnarray}
Since this integral cannot be done in closed form unless $c=-1$ (the limiting case of a sphere), one must numerically integrate this. Introduce the reciprocal of the angle $\theta$ measured in cycles (i.e., the reciprocal of $\theta/(2\pi)$)
\beq
  N(\chi,\beta_0) 
= \frac{2\pi}{G(\chi,\beta_0)} \,,\quad \beta_0 < \beta_{\rm(crit)}
\,.
\eeq
This can be used to numerically determine the values of the initial angle $\beta_0$ which correspond to closed orbits.

\begin{figure}[t] 
\typeout{*** EPS figure torus bound periodic geo [m,n;1] grids}
\begin{center}
\includegraphics[scale=0.3]{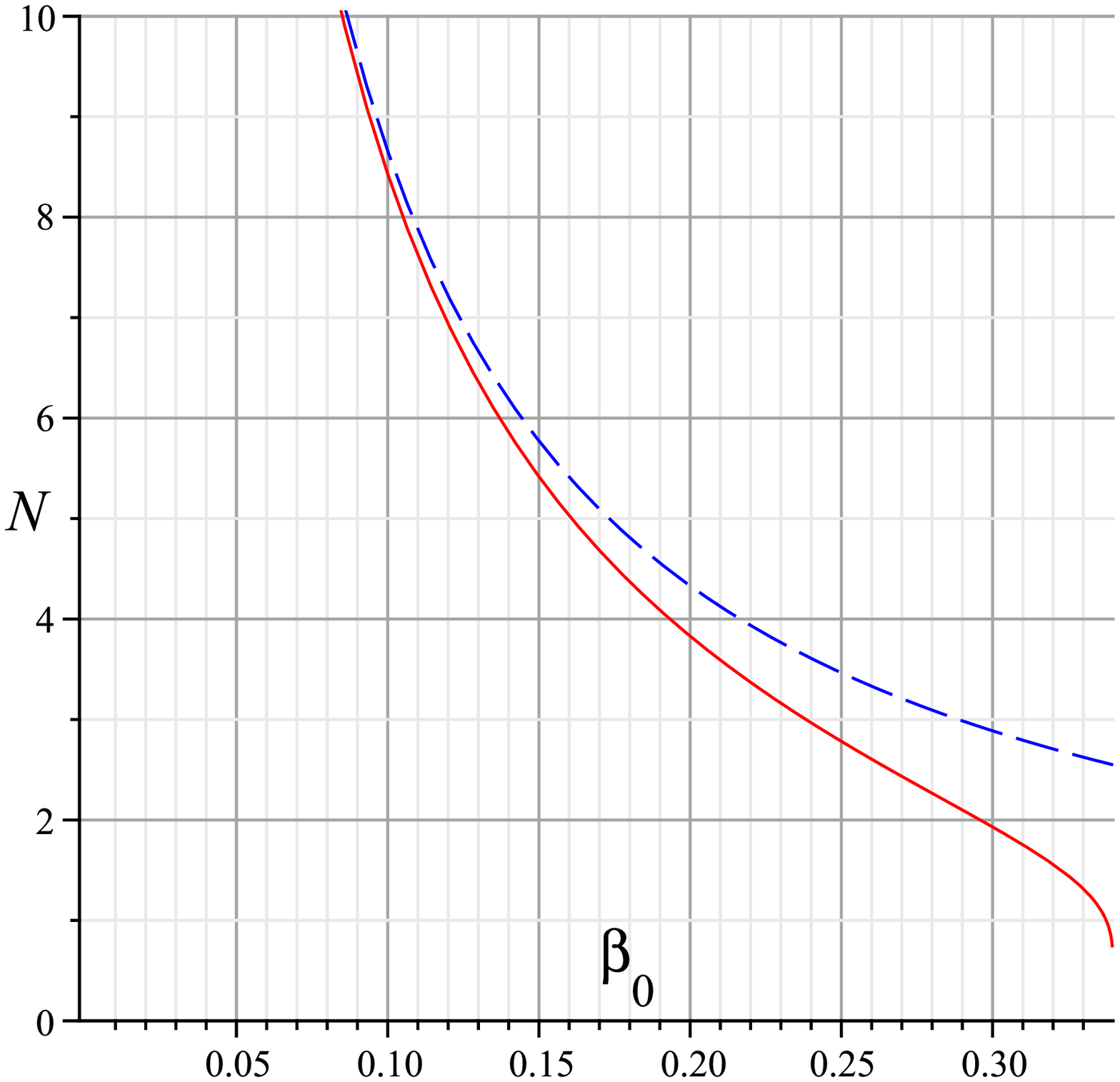}
\includegraphics[scale=0.3]{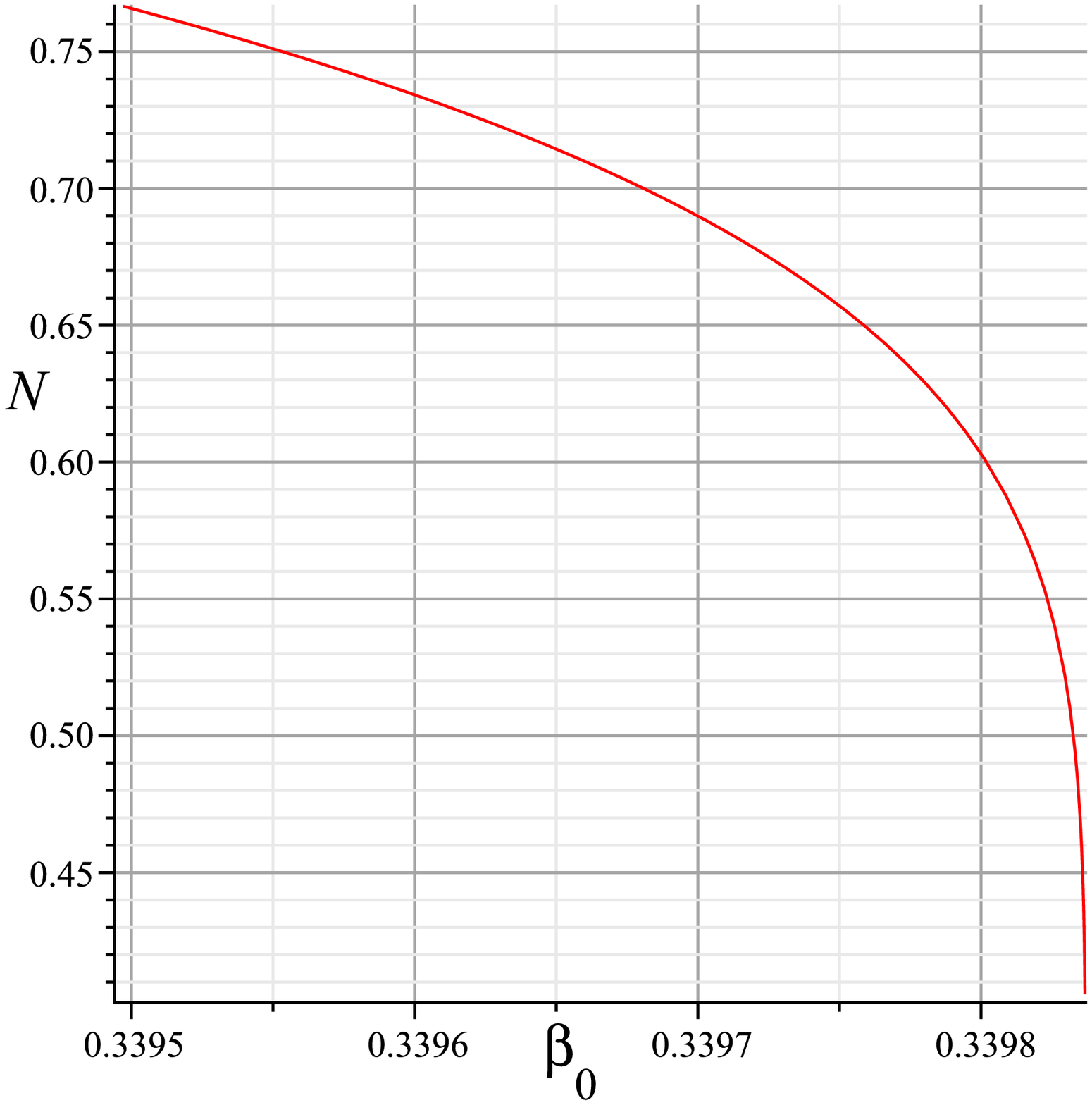}
\end{center}
\caption{
The number $N(2\pi,\beta_0)$ of radial oscillations per azimuthal revolution for the unbound geodesics on the unit ring torus plotted versus the initial angle $\beta_0$ first over the interval $[0,\beta_{\rm(crit)}]$ and then for $[0.9999 \beta_{\rm(crit)},\beta_{\rm(crit)}]$, where $\beta_{\rm(crit)}\approx 0.3398$. The dashed line shows the asymptotic $k/\beta_0$ approximation for small angles.
The rational values $N(2\pi,\beta_0)=m/n$ correspond to the unbound closed geodesics of type $[m,n;1]$.
}
\label{fig:unboundgrid}
\end{figure}

For an unbound geodesic the value of this quantity after one radial oscillation is the theta frequency function $N(2\pi,\beta_0)=2\pi/G(2\pi,\beta_0)$ associated with the theta period $G(2\pi,\beta_0)$. This quantity is the number of times the azimuthal increment occurring during that complete oscillation fits into a full azimuthal revolution, i.e, the number of radial oscillations per azimuthal revolution. If this is a rational number $m/n$
\beq
  N(2\pi,\beta_0) =\frac{m}{n}\,,
\eeq
then after $m$ radial oscillations, $n$ azimuthal revolutions will take place, corresponding to an unbound geodesic of type $[m,n;1]$. Fig.~\ref{fig:unboundgrid} shows a plot of $ N(2\pi,\beta_0) $
over the entire interval $\beta\in[0,\beta_{\rm(crit)}) $ of values corresponding to the unbounded geodesics moving up initially in the positive $\theta$ direction. The graph rises very quickly due to a vertical asymptote at $\beta_0=0$, while it decreases extremely slowly to its limiting value $\lim_{\beta_0\to\beta_{\rm(crit)}} N(2\pi,\beta_0)=0$ where it approaches a vertical tangent, as shown in a closeup of the endpoint behavior.

The integral $F(\pi,\beta_{\rm(crit)})$ is an improper integral whose integrand denominator goes to zero at the right endpoint, and it has an infinite value corresponding to the infinite number of azimuthal revolutions which occur as the critical geodesic asymptotically approaches the inner equator. The reciprocal $N(\pi,\beta_0)$ therefore goes to zero as $\beta_0\to\beta_{\rm(crit)}$ from the left, but very slowly.
To get an idea of how slowly,
set $\beta_0=\beta_{\rm(crit)}$ (which satisfies $(c+1)\sin\beta_{\rm(crit)}=c$),
and evaluate the Taylor expansion about $\chi=\pi$ of the reciprocal of the integrand of $F(\pi,\beta_{\rm(crit)})/(2\pi)$ to find that the leading behavior of the integrand at its endpoint $\chi=\pi$ is
\beq
   \frac{1}{2\pi c^{1/2} (\pi-\chi)}
\eeq
so its antiderivative behaves like
\beq
   \frac{1}{2\pi c^{1/2}} \ln[(\pi-\chi)^{-1}]\,.
\eeq
Thus if in the discrete numerical plotting of this function at the endpoint, if the last sampled point is $\pi-\chi=10^{-10}$, say, then the value of the antiderivative there for the unit ring torus ($c=1$) is only about
\beq
   \frac{1}{2\pi} \ln[(\pi-\chi)^{-1}]\approx 3.3
\,.
\eeq 
The reciprocal of this gives a rough idea of the value of $N$ reached by this last sampled point. This means numerically finding rational values of $N$ for even relatively large proper fractions where $\beta_0$ is close to the critical value will be difficult.
For example, for the unit ring torus where $\beta_{\rm(crit)}\approx 0.3398369094$ ($19.471^\circ$),
one easily finds $\beta_0=0.2382795502$ ($13.7^\circ$) for the $[3,1;1]$ geodesics
and $\beta_0=0.0.3226432999$ ($18.5^\circ$) for the $[3,2;1]$ geodesics, 
but the result $\beta_0=0.0.3226432999$ ($19.454^\circ$) for the $[3,4;1]$ geodesic takes a considerable time, while 
the value $\beta_0=0.3395532232$ ($19.469^\circ$ for the $[3,5;1]$ geodesic takes too long for a simple root finder to obtain and one must resort to trial and error.

\begin{figure}[t] 
\typeout{*** EPS figure torus bound periodic geo bound grids}
\begin{center}
\includegraphics[scale=0.5]{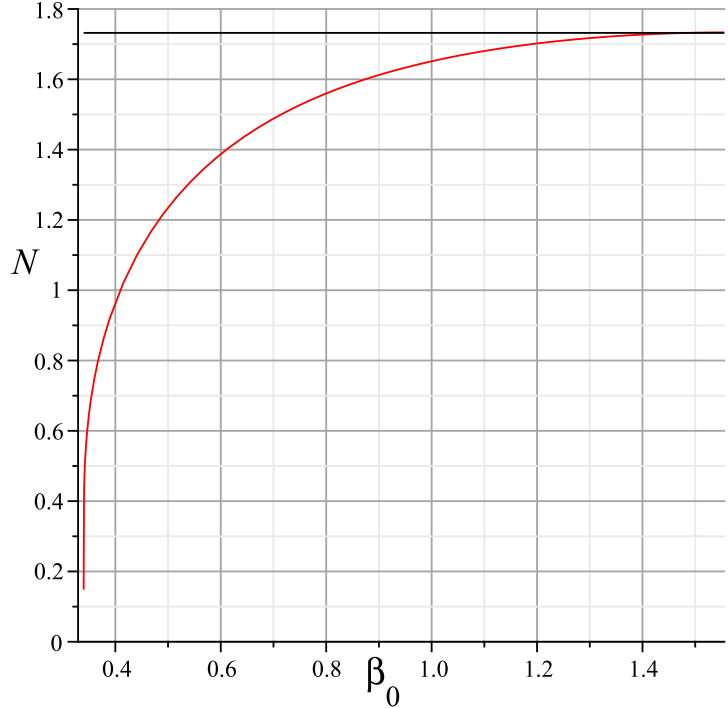}
\end{center}
\caption{
The number $N(2\pi,\beta_0)$ of radial oscillations in one azimuthal revolution versus the initial angle $\beta_0\in (\beta_{\rm(crit)},\pi/2]$, here shown for the unit ring torus with the upper bound $N(2\pi,\pi/2)=(c+2)^{1/2}=\sqrt{3}\approx 1.73$. The rational values $N(2\pi,\beta_0)=m/n$ correspond to the bound closed geodesics of type $[m,n;0]$ with $m$ radial oscillations during $n$ azimuthal revolutions.
}
\label{fig:boundgrid}
\end{figure}

The azimuthal integral $\theta=G(\chi,\beta_0)$ for the orbits of the unbound geodesics can be easily approximated in the limit of very small initial angles where a small increment of the azimuthal angle occurs during one radial oscillation. Under the condition $\ell/\sqrt{2E}=(a+b)\sin\beta_0 \ll1$, the integral approaches
\beq
\theta = \int_0^\chi \frac{\ell}{\sqrt{2E}} \frac{b\,d\chi}{R(b\chi)^2} \,,
\eeq
which has a complicated exact value but the increment over one revolution $-\pi\le \chi\le \pi$ is relatively simple (computer algebra system exercise)
\beq
\Delta \theta = \frac{2\pi ab(a+b)}{(a^2-b^2)^{3/2}} \beta_0 \,.
\eeq
This leads to
\beq
\beta_0\to0:\qquad
  N(2\pi,\beta_0)\to  
\frac{(a^2-b^2)^{3/2}}{ab(a+b)} \frac{1}{\beta_0} 
= \frac{c^{3/2}(c+2)^{1/2}}{(c+1)} \frac{1}{\beta_0}
\,,
\eeq
which nicely fits the approach to the vertical asymptote for small angles as shown in Fig.~\ref{fig:unboundgrid}.

For a bound closed geodesic, the increment $\Delta\theta$ in $\theta$ during one radial oscillation cycle is four times the increment from $\chi=0$ to $\chi=\chi_{\rm max}$, so one must modify the number of radial oscillations per azimuthal revolution function to
\beq
  N_{\rm(bound)}(\beta_0)  = \frac{2\pi}{4G(\chi_{\rm max}(\beta_0),\beta_0)} 
\le (c+2)^{1/2}
 \,,\quad \beta_{\rm(crit)}<\beta_0 <\frac{\pi}{2} 
\,.
\eeq
This upper bound for the theta frequency of the bound geodesics comes from the limiting case of small radial oscillations about the outer equator which is the largest radius circle of revolution in the torus. Discussed more at length in Appendix C,
as the amplitude of the oscillation increases (increasing $\chi_{\rm(max)}$), the conservation of angular momentum at the smaller radii of revolution increases the azimuthal angular velocity compared to the radial motion, thus increasing the ratio of azimuthal revolutions per radial oscillation. Fig.~\ref{fig:spray} of Appendix B shows this behavior for the unit ring torus.

The rational values $N_{\rm(bound)}(\beta_0) =m/n$ correspond to the bound closed geodesics of type $[m,n;0]$ with $m$ radial oscillations during $n$ azimuthal revolutions. This condition must be solved numerically for a given ratio $m/n$ 
bounded above by the condition
$$
m/n <(c+2)^{1/2}\,.
$$
For the unit ring torus where $(c+1)^{1/2}=\sqrt{3}\approx1.73$, for example, the (inequivalent) bound closed geodesics can be organized by increasing $n$ and then $m$
\beq
 \frac11;
 \frac12,\frac32;
 \frac13,\frac43,\frac53;
 \frac14,\frac34,\frac54;
 \frac15,\frac25,\frac35,\frac45,\frac65,\frac75,\frac85; \ldots \,.
\eeq
For any bound geodesic $r$ is implicitly a periodic function of $\theta$ with period $4G(\chi_{\rm max}(\beta_0),\beta_0)$. For the closed such geodesics, the above condition simply requires that this period be a rational multiple of $2\pi$.

\begin{figure}[t] 
\typeout{*** EPS figure torus bound periodic geo [1,1;0], [1,2;0]}
\vglue-1.5cm
\begin{center}
\includegraphics[scale=0.3]{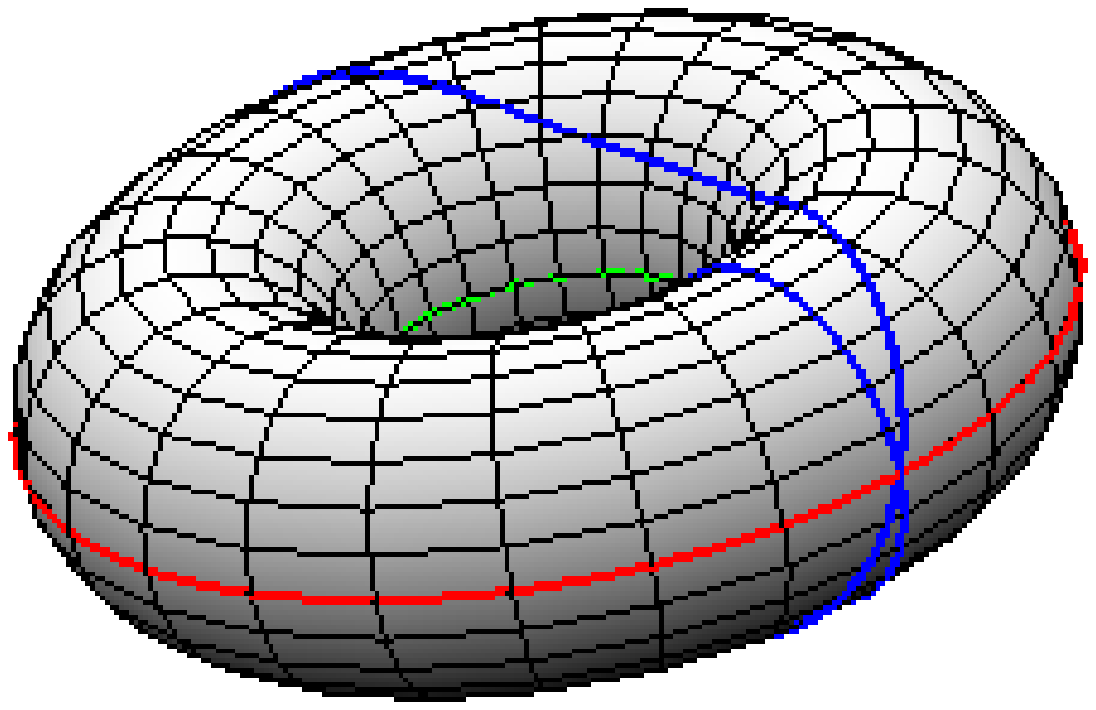}
\includegraphics[scale=0.3]{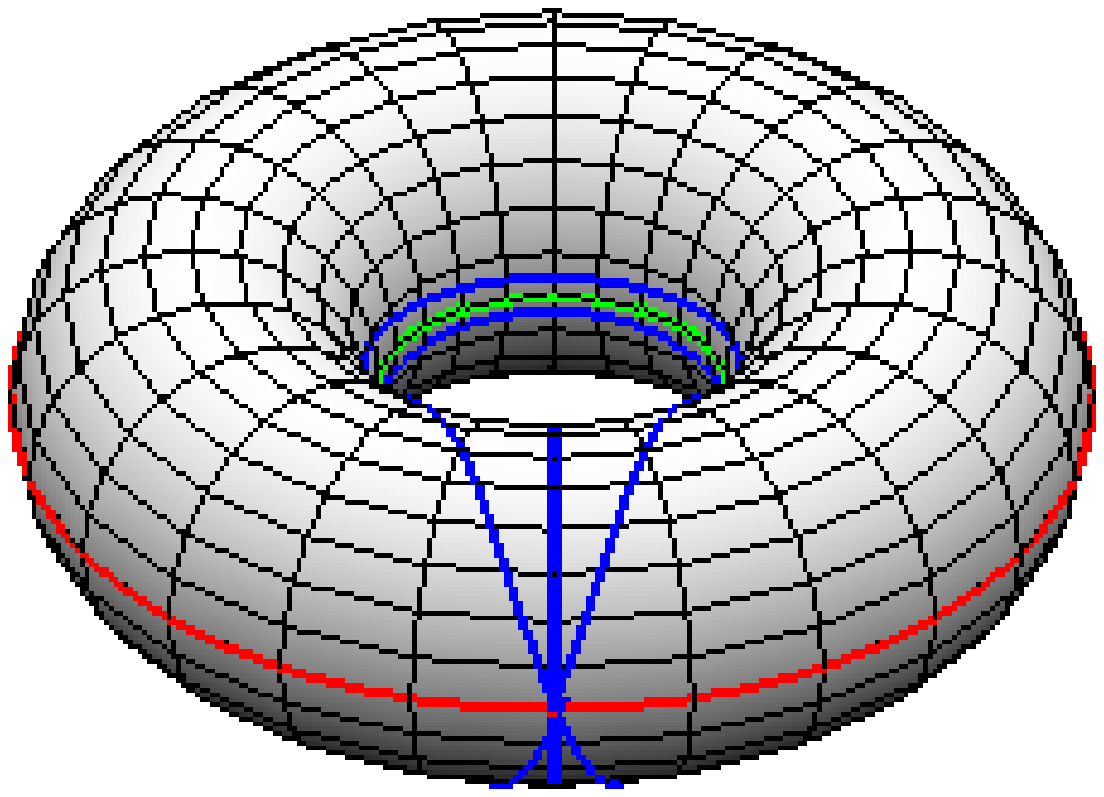}
\end{center}
\vglue-1.5cm
\caption{
Bound periodic orbits $[1,1;0]$, $[1,2;0]$  for the unit ring torus corresponding to initial angles of $\beta_0\approx 0.4097, 0.3422 \,(\approx 23.5^\circ, 19.6^\circ)$ and maximum radius $\chi_{\rm max}\approx 143.6^\circ, 173.4^\circ $ and lengths $15.3, 21.9 $.
Orbits with initial angle $\beta_0\le41.8^\circ$ reach the North Polar Circle,
while orbits with initial angle $\beta_0\le19.47^\circ$ are unbound, passing the inner equator if the inequality is satisfied.
}
\label{fig:torusperiodicgeos110-120}
\end{figure}

\begin{figure}[t] 
\typeout{*** EPS figure torus bound almost periodic geo [7,5;0]}
\begin{center}
\includegraphics[scale=0.5]{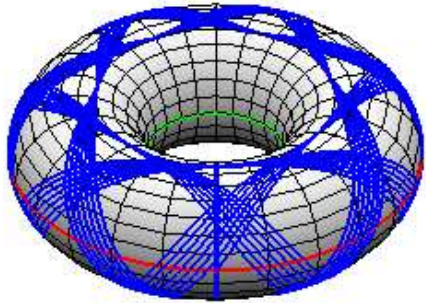}
\end{center}\vglue-3cm
\caption{
An almost closed bound geodesic close to the $[7,5;0]$ closed geodesic ($\beta_0=0.6124$, 35.08 degrees) with a slightly higher initial angle ($\beta_0=0.614$, 35.18 degrees) that causes the geodesic to return to the outer equator near the initial meridian a bit too soon, and therefore causing the nodal point of return near the initial meridian to precess in the opposite direction of the orbit (the orbit increases the azimuthal angle, the precession decreases the azimuthal angle of the node). A slightly smaller initial angle would lead to precession in the same direction as the orbit. Trial and error adjustment of the (slightly smaller) approximate initial angle from the integral condition enables one to narrow in on the correct value to as many digits as one has patience for, subject to the unknown accuracy of the numerical integration scheme of course.
}
\label{fig:torusperiodicgeo75precession}
\end{figure}

Figs.~\ref{fig:torusperiodicgeos32} and \ref{fig:torusperiodicgeos110-120}  show the first three geodesics of this list.
It is relatively easy to find the initial angles for the larger permitted ratios $m/n$ for sufficiently small integers, but as the ratio $m/n$ gets smaller, the limiting vertical tangent in the graph of $N_{\rm(bound)}(\beta_0)$ as $\beta_0\to\beta_{\rm(crit)}$ from the right creates problems for accurate determination of the roots where $N_{\rm(bound)}$ is a small rational number. Because its reciprocal is an improper integral with a positive integrand that goes to infinity at the critical angle, numerical evaluation will underestimate the result of the integral (because numerically it will only sense some highest value of the integrand in the discrete approximation and miss the infinite contribution). Thus the true angle which gives a particular rational value to $N_{\rm(bound)}(\beta_0)$ will have to be a bit bigger to compensate for the lower value of the approximate integral. Thus one can take this underestimate as a starting angle in the full numerical geodesic solver and increase it slightly by trial and error until one closes the orbit. Even for larger values of the ratio $m/n$ where one can accurately determine the corresponding initial angle, if $m$ or $n$ gets too large, one is of course limited again by numerical error in the full geodesic solver since too many radial oscillations or azimuthal revolutions translates into accumulated numerical error, so the orbit may not appear to close. When it fails to close by just a little, one sees a constant rate of ``precession" of the node (where the orbit crosses the outer equator $r=0$) since the crossing angle at the outer equator is the same each time it crosses nearby and hence it repeats the original trajectory but translated around the torus slightly by the amount by which the return position to the outer equator misses the original meridian of its departure. This is illustrated in 
Fig.~\ref{fig:torusperiodicgeo75precession}.

\section{Lengths of geodesics}

Since increments of the affine parameter are only proportional to the arc length of a geodesic segment, one must evaluate the arc length from the affine parameter increment using the constant of proportionality $(2E)^{1/2}$, which is a function of the
energy or equivalently of the initial angle. If one has only incomplete information not including the affine parameter increment, like the initial angle for a closed geodesic,  the arc length integral remains to be evaluated, say as a function of the interval of radial coordinate values describing the geodesic segment. 
The corresponding affine parameter integral given in Eq.~(\ref{eq:lambdaofr}) can be evaluated exactly in term of elliptic functions, but the resulting complicated formula is not very enlightening. With a slight modification, this formula may be used to
evaluate the arc length of one complete period of the motion for a closed geodesic.

To evaluate the arc length of one period of one of the unbound closed geodesics of the class $[m,n;1]$, one can first evaluate the arc length of one radial loop ($0\le r \le 2\pi b$ or $0\le \chi \le 2\pi$) and then multiply by the number $m$ of radial loops. The arc length of one radial loop follows from integrating 
\beq
  \frac{ds}{dr} =   \frac{ds/d\lambda}{dr/d\lambda}  = \frac{(2E)^{1/2}}{dr/d\lambda}\,,
\eeq
namely
\begin{eqnarray}
  L &=& \int_0^{2\pi b} \frac{(2E)^{1/2} \,dr}{(2E-\ell^2/R^2)^{1/2}}
   = b \int_0^{2\pi } \frac{R \,d\chi}{(R^2-\ell^2/(2E))^{1/2}} 
\nonumber\\
 &=&    b \int_0^{2\pi } \frac{R \,d\chi}{(R^2-R(0)^2 \sin^2\beta_0)^{1/2}} 
\,.
\end{eqnarray}
For example for the 7 loop geodesic with energy $E=3.942$ on the unit ring torus, the length of one loop is 6.45 and of the entire 7 loops: 45.12. This compares to the circumferences of the inner equator, the Polar Circles and the outer equator: $[2\pi(1),2\pi(2),2\pi(3)]\approx [6.28,12.57,18.85]$. These integrals are easily done numerically with a computer algebra system. Notice that for very high energy (very small initial angle $\beta_0\to0$), this reduces to the radial circumference of the torus.

For any bound geodesic, one can easily evaluate the length of one transit departing from the initial point on the outer equator, traveling to a given turning point, returning to cross the outer equator and going to the other turning point and back again by taking 4 times the length from the outer equator to one turning point
\beq
 L   = 4\, b \int_0^{\chi_{\rm(max)}}  \frac{R \,d\chi}{(R^2-R(0)^2 \sin^2\beta_0)^{1/2}} \,,
\eeq
where
\beq
   \chi_{\rm(max)}
                = \arccos\left( (2+c)|\sin\beta_0|-(1+c) \right)
\,.
\eeq
For example the [1,1;0] closed geodesic shown in Fig.~\ref{fig:torusperiodicgeos110-120} on the unit ring torus has length 15.26. 

One can also categorize the closed geodesics by the number of self-intersec\-tions they exhibit, as initiated by Irons \cite{irons}. To examine geodesic crossings without trying to follow a geodesic around the 3-$d$ torus in rotatable computer graphics images, it is useful to map the torus onto the unit square by introducing two angles measured in cycles rather than radians
\beq
   (\Xi,\Theta)=\left( \frac{\chi}{2\pi},\frac{\theta}{2\pi} \right) \,, \ \mbox{with}\ 
  0\le  \chi,\theta\le 2\pi \longrightarrow
  0\le \Xi,\Theta \le 1\,.
\eeq

For comparison the grid of all positive pairs $[m,n]$ describing the distinct closed geo\-desics of the flat torus is briefly discussed in the Appendix. Like the flat torus, the $[m,n;1]$ closed geodesics on a ring torus exist for all positive pairs $(m,n)$ and have no self-intersections, but the $[m,n;0]$ closed geodesics are constrained by the theta frequency condition $m/n<(c+2)^{1/2}$ of the limiting geodesics near the outer equator, which have the largest period  of the bound geodesics, and self-intersections are allowed.
As one gets closer and closer to the critical geodesic which asymptotically approaches the inner equator, the closeness to the unstable equilibrium there allows an arbitrarily large number of azimuthal revolutions near it before the geodesic flies off towards the outer equator again, which can be arranged to happen to hit any point on the outer equator. Thus one can get all pairs $(m,n)$ for $m/n\to0$.
Numerically this is problematic, however, since numerical error is hard to control near an unstable equilibrium. 

Given that only the bound geodesics can have self-intersections, comparison of the low integer geodesics of type $[m,n;0]$ reveals the following apparent crossing rules as a function of the number $m$ of radial oscillations. For the bound geodesics on the unit ring torus starting at the origin of coordinates with initial angle $\beta \in (\beta_{\rm(crit)},\pi/2)$ orbiting in the increasing $\theta$ direction, looking down from above at the upper hemitorus one finds by inspection of many cases either 0, 1, or 2 equally spaced series of self-crossings at 0, 1, or 2 distinct values $\chi_{\times}$
of the radial angle $\chi \in[0,\chi_{\rm(crit)}]$. The equal spacing increment is $2\pi/m$.

\vskip\baselineskip

$n=1$: No self-crossing point radii.

$n=3$: 1 self-crossing point radius: $\chi_{\times}$. The series is shifted right from  $\theta=0$  by $(2\pi/m)/4$.

$n>3$ odd: 2 self-crossing point radii: $\chi_{\times 1}<\chi_{\times 2}$. The  2 series are shifted right/left from $\theta=0$ by  $(2\pi/m)/4$.

$n=2$: 1 self-crossing point radius $\chi_{\times}=0$.

$n>2$ even: 2 self-crossing point radii:   first series at $\chi_{\times 1}=0$, second series $\chi_{\times 2}$
shifted right from $\theta=0$ by  $(2\pi/m)/2$.

\vskip\baselineskip

\noindent
It would be nice to understand how these rules come about, and to evaluate the angles $\chi_\times$ at which these crossings occur.

Finally while it is convenient to describe a geodesic as the curve between two points which minimizes the distance between them on the surface, it is also important to recognize that the converse is not true. While a geodesic connecting two points is a critical point of the arc length function on the space of paths between the two specified points, as follows from the variational approach to this problem, there can be many geodesics of differing lengths connecting any two points on a surface. In the flat plane, the geodesic is the unique straight line between two points, but on a sphere although there is a unique great circle containing any two non-antipodal points, there is still a choice of the longer or shorter segment of the great circle to follow from one point to the other, and for antipodal points there are an infinite number of great circles to choose from. For the torus the situation is even more extreme: through any two points on the torus, there are an infinite number of both bound and unbound geodesics which connect them. However, for most point pairs there is still a minimal length geodesic connecting them. While one can try to aim a geodesic departing from the first point to hit the target point using trial and error, it is interesting that completely different approaches can be taken to numerically converge on this 
connecting path \cite{baek,murat,lopez}.

An interesting illustration of the minimal geodesic question regards the minimal length geodesic which connects two distinct points on the outer equator. According to the approximate analysis above, geodesics departing from $(r,\theta)=(0,0)$ with an initial angle differing very little from $\beta=\pm \pi/2$ oscillate about the outer equator with half period $\Delta\theta=\pi/(c+2)^{1/2}$, thus returning to the outer equator first at $\theta=\Delta\theta$. Examining numerical geodesics for such initial conditions then shows that for the actual geodesics, as the amplitude of the radial oscillation increases, the actual half period increases slowly, so that the return point migrates slowly beyond $\theta=\Delta\theta$, but the length of this geodesic segment remains smaller than the length of the direct route between the endpoints along the outer equator. Thus for two points on the outer equator with separation equal to or less than the half period $\Delta\theta$, the outer equator segment between them is the shortest geodesic between them, but for points with separation greater than this critical length, the shortest geodesic is instead a half radial oscillation connecting the two points. Antipodal points on the outer equator are the extreme case of this type, where the $[1,1;0]$ geodesic which connects those points is significantly shorter than half the circumference of the outer equator.

The horizontal ``focusing" of the bound geodesics which emanate from the origin of coordinates on the outer equator, like the crossing of great circles on the sphere at antipodal points, is a consequence of the positive curvature of the torus on the outer ring between the Northern and Southern Polar Circles. On the inner ring between these Polar Circles, the curvature is negative and geodesics from the outer equator which reach this region instead diverge from one another before crossing on the inner equator due to the topology of the torus. These geodesics begin crossing each other at the critical length $L_r=\pi b$ equal to half the radial circumference of the torus, compared to the azimuthal critical length $L_\theta=(a+b)\Delta\theta =\pi (a+b)/(c+2)^{1/2}$, leading to the ratio $L_\theta/L_r=(c+2)^{1/2}$. Irons \cite{irons} has evaluated the various curvature quantities on the torus, but curvature is a long story that cannot be told here, but which is described in Appendix C.

\section{Kepler problem?}

The beauty of mathematics is that many seemingly different problems actually share the same mathematical foundation and hence can be treated with the same methods to understand them. From geodesics on these special surfaces of revolution,  by relaxing the condition $g_{rr}=1$ to $g_{rr}=f(r)$ we can easily consider general surfaces of revolution for which the radial coordinate does not measure arc length \cite{gray3}. This enables us to handle the simple parabolic surface of revolution $z=x^2+y^2=r^2$, for example, where $g_{rr}=1+4r^2$. However, we can also handle the simplest example of a surface of revolution, the flat plane $z=0$, but in the framework of polar coordinates, where resolving the geodesic equations to reproduce straight lines is complicated but doable. Polar coordinates are simply not adapted to the geodesic structure of the flat plane. The real payoff, however, is to consider the simple addition of an external rotationally symmetric force to these same flat plane equations, enabling us to handle the nonrelativistic Kepler problem with the same approach, and indeed the closely related relativistic Kepler problem of planar motion in a spherically symmetric vacuum gravitational field arising either from a spherical body or a black hole. In both cases, the polar coordinates are essential to resolving the equations of motion.

\begin{figure}[t] 
\typeout{*** EPS figure Kepler potential}
\begin{center}
\includegraphics[scale=0.35]{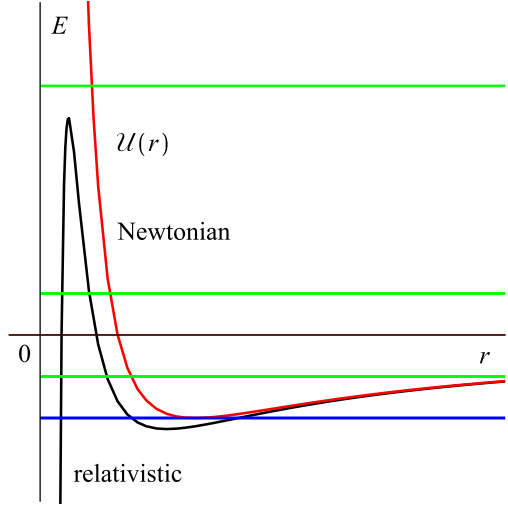}
\end{center}
\caption{
A typical Kepler potential for the nonrelativistic (Newtonian gravity) and relativistic (general relativity) cases, showing the ``centrifugal barrier" near the source (small $r$). Circular orbits occur at the extrema of the potential. Bound orbits have  energy $E<0$, while unbound orbits have $E\ge0$. The relativistic potential turns over close to the source where stronger gravitational attraction overcomes the centrifugal force. The lowest two energy levels shown correspond to bound orbits with a minimum and maximum radius, while the upper two are unbound, but in the relativistic case, the highest energy level orbit overcomes the centrifugal barrier and is captured by the source, while the next lower energy orbit is reflected from the centrifugal barrier at the minimum radius in both cases (the outer portion for the relativistic case); the relativistic case also contains ``trapped" orbits inside the barrier. 
} 
\label{fig:keplerpotential}
\end{figure}

For the flat plane case $R(r)=r$ and the energy equation becomes
\beq
 E = \frac12 \left(\od{r}{\lambda}\right)^2 
    + \frac12 \frac{\ell^2}{r^2} \,.
\eeq
To incorporate an actual rotationally symmetric radial force $F_r(r)=-d U_r(r)/dr$ into the mix, we just need to add its potential $U_r(r)$ to the effective potential term already present
\beq
 E = \frac12 \left(\od{r}{\lambda}\right)^2 
    + \frac12 \frac{\ell^2}{r^2}+U_r(r) \,.
\eeq
The new total potential then governs the problem in the same way that the effective potential alone did the job in our previous discussion. With a potential of the form
\beq
  U_r(r) = -\frac{k_1}{r}-\frac{k_2}{r^3} \rightarrow F_r(r)= -\frac{k_1}{r^2}-\frac{3k_2}{r^4} \,,
\eeq
we can handle the nonrelativistic Kepler problem for motion in an inverse square force if $k_2=0$ (attractive like gravitation if $k_1>0$) or the relativistic Kepler problem ($k_1>0,k_2>0$, see \cite{MTW}), see Fig.~\ref{fig:keplerpotential}. The additional inverse cubic potential term in the latter case corresponds to an additional inverse quartic force that is responsible for the precession of the orbits which in its absence are conic sections.

Amusingly enough, the radial oscillations about the outer equator of the torus are directly mirrored in the Kepler problem, and the frequency one calculates for the oscillations about a stable circular orbit is called the epicyclic frequency, an important astrophysical quantity.
While these more general problems are as easily dealt with as the torus geodesics, they are perhaps best left to a more comprehensive discussion elsewhere  \cite{drbob}. However, a slight adjustment of the Maple worksheet for the torus allows experimentation with these kinds of orbits as well. With a little more work, one can include a light pressure drag force to model the effect of stellar radiation pressure on orbiting dust particles, the so called Poynting-Robinson effect  \cite{pre}.

Perhaps the most important concept highlighted here that is missing in most undergraduate mathematics programs is energy. There is usually just no time in an introductory course on differential equations to introduce it (and the standard textbooks do not mention it), so whatever vague idea mathematics students have of it from high school or possibly from college freshman physics is never connected concretely to differential equations and the natural way it arises there. Even in the present article, we have danced around this connection without directly confronting it either, but hopefully its utility at least in this problem is clear. The usual approach to geodesics through a variational approach is the natural way to make this connection \cite{feldman}. This is an interesting advanced topic either for an undergraduate mathematics or physics setting. For advanced undergraduate mechanics in the physics setting, this approach to constrained motion (motion on a surface) is a natural gateway for understanding the tools of general relativity, as developed by Walecka \cite{walecka}, for example.

\section{What next?}

It would be sad to shine the light on an interesting problem and completely exhaust its potential for entertainment. 
There is much left to explore. On the mathematical side,
the spectrum of closed geodesics on the torus depends on the shape parameter. Passing to the horn torus $c=0$ from the ring tori, the local maximum in the periodic potential at the inner equator moves off to infinity, and the unbound nonradial geodesics disappear. All the nonradial geodesics are repelled from the single point inner equator by this angular momentum barrier.
For two-point self-intersecting tori $-1< c< 0$, this infinite angular momentum barrier splits into two barriers surrounding the two self-intersection points on the symmetry axis, leading to two disjoint families of nonradial bound geodesics, those confined to the outer apple surface and those confined to the inner lemon surface. 
%
%
For the shape parameter $c=2$ (for example: $(a,b)=(3,1)$ so the ratio of the outer equator radius to inner equator radius  is 2), exactly two radial oscillations about the outer equator fit into one azimuthal revolution in the limit of zero amplitude oscillations. This means that for nonzero amplitudes where the exact frequency will be slightly different from 2, the nodal points where the oscillations intersect the outer equator will slowly migrate or ``precess" in the physics language. What is the precession frequency? This can be analyzed by keeping the next term beyond the quadratic one in the Taylor expansion of the potential about its minimum. (Similarly one can try to analyze the precession of the node for near miss closed geodesics to calculate the incrementing angle of precession described above.)
For the pinched torus, all radial geodesics pass through the pinching point at the origin, but all of them have a vertical tangent vector, so one looses uniqueness of the initial value problem for the geodesics there; one can match up radial geodesics of different azimuthal angle on either side of the pinching point. What does this do to the torsion and differentiability of those curves?

On the physical side, one can easily explore motion in various classical mechanical force potentials as hinted at for the Kepler problems. In an advanced mechanics course, where only formula manipulation was possible when I was a student, there are exciting possibilities for student experimentation.

Finally the obvious must be pointed out. None of this would be possible without the incredible tool now at the disposal of faculty and students alike: the computer algebra system on our personal computers. Both of the two leading such software packages Maple and Mathematica offer to motivated students an enormous potential for exploring mathematical and physical problems, if given a little guidance from faculty. The problem for undergraduate teaching is to find the opportunities to engage students in picking up this tool in parallel with their gradual increase in theoretical sophistication in either a mathematics or physics setting.

\appendix
\section{The flat torus closed geodesics}

\begin{figure}[t] 
\typeout{*** EPS figure flat torus 23}
\begin{center}
\includegraphics[scale=0.35]{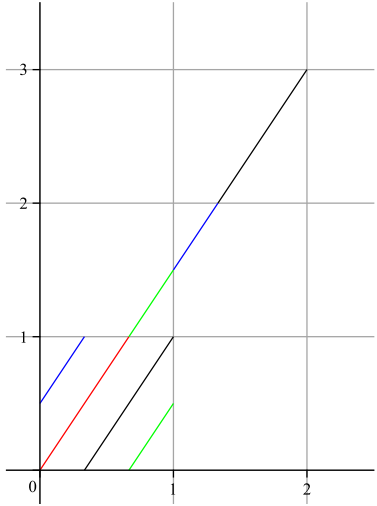}
\end{center}
\caption{
The $[2,3]$ closed geodesic through the origin before and after imposing the periodic boundary conditions for the flat torus on the unit square.
} 
\label{fig:torusflat23}
\end{figure}

\begin{figure} 
\typeout{*** EPS figure flat torus up to 6}
\vglue-0.5in
\begin{center}
\includegraphics[scale=0.35]{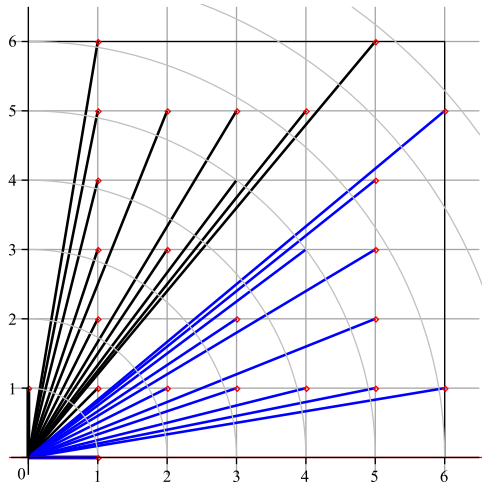}
\end{center}
\vglue-0.3in
\caption{
The closed geodesics on the flat torus on the unit square characterized by the period pairs  $[m, n]$ for $ m, n < 7$ before imposing the periodic boundary conditions.
} 
\label{fig:torusflat66}
\end{figure}

For comparison it is useful to consider the spectrum of closed geodesics for the simpler flat torus. The simplest realization of the flat torus is the plane with points identified  if their Cartesian coordinates differ by integers, thus reducing distinct points to the unit square in the first quadrant with opposite boundary line segments identified.
One can classify these geodesics by the pair of periods $[m,n]$ of the number of revolutions in the horizontal and vertical directions respectively. Fig.~\ref{fig:torusflat23} shows a typical geodesic of type $[2,3]$ through the origin returning the origin at the point $(2,3)$ in the plane, with length $\sqrt{2^2 +3^2}=\sqrt{13}$.
Fig.~\ref{fig:torusflat66} shows in the same way a representation of the lattice of all the distinct closed geodesics through the origin of type $[m,n]$ for $m,n<7$, with length $\sqrt{m^2+n^2}$.

\begin{figure} 
\vglue-0.3in
\typeout{*** EPS figure ftorus geo rays in tangent space}
\begin{center}
\includegraphics[scale=0.4]{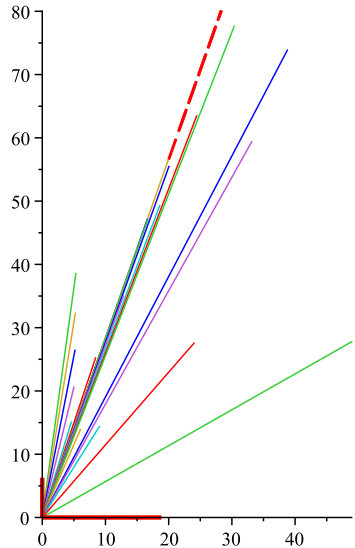}
\end{center}
\vglue-0.3in
\caption{
The closed geodesics on the unit ring torus with initial angles $0\le \beta_0 \le \pi/2$ characterized by the triplets  $[m,n;p]$ for $ m, n < 7$, $p=0,1$, represented in the tangent plane at the outer equator by a line segment of length equal to the arc length of the geodesic and with the same angle $\beta_0$ measured clockwise from the positive vertical axis. Included are the vertical meridian and the horizontal outer equator of lengths $2\pi$ and $6\pi$ respectively. The unbound (top) and bound (bottom) geodesics are separated by the critical angle $\beta_0=0.3398369094$ (19.5 degrees) indicated by the dashed line.
} 
\label{fig:torusgeorays}
\end{figure}

One can reinterpret Fig.~\ref{fig:torusflat66} as a representation of the actual closed geodesics in the tangent space at the origin using the exponential map. This map associates a point in the tangent space at a point $P$ with a point $Q$ on a geodesic through it by assigning to it a tangent vector with the same direction as the tangent vector to the geodesic at the point $P$, but with a length equal to the arc length from $P$ to $Q$ along that geodesic. One can therefore reproduce this figure for the bound and unbound geodesics on the curved unit ring torus. One only needs to assemble the pairs $(\beta_0,L(\beta_0))$ of initial angles and corresponding lengths for the low integer  $(m,n)$ pair bound and unbound closed geodesics and plot them in the vertical tangent plane at the origin of coordinates. Fig.~\ref{fig:torusgeorays} shows this result.
One can read about the exponential map in most mathematical textbooks which discuss geodesics \cite{exp}.

\section{Variational approach to the geodesic equations}

One can side step all of the Riemannian geometry discussion of covariant differentiation and parallel transport by using the variational approach which produces the geodesic equations by extremizing an action functional on the space of all curves connecting any two fixed points on the torus. One can either extremize the arc length of the curve, which is the integral of the length of the tangent vector in any parametrization of the curve, namely the speed function
\beq
 Action_1 = \int_c ds = \int_c \frac{ds}{d\lambda} \,d\lambda 
= \int_c \sqrt{\left(\frac{dr}{d\lambda}\right)^2+ \left(\frac{d\theta}{d\lambda}\right)^2} \,d\lambda \,,
\eeq
which is clearly independent of a change of parametrization,
or extremize the integral of half the length squared of the tangent vector (the factor of two is for convenience)
\beq
 Action_2 = \frac12 \int_c \left(\frac{ds}{d\lambda}\right)^2 \,d\lambda 
= \frac12 \int_c \left(\frac{dr}{d\lambda}\right)^2+ R^2\left(\frac{d\theta}{d\lambda}\right)^2 \,d\lambda \,,
\eeq
which is equivalent to the previous case only for affinely parametrized curves for which the speed $ds/d\lambda$ is constant. These facts require a separate study of the topic called the calculus of variations \cite{calcvariations}. 
The integrand is called a Lagrangian function, and it is a function of the curve and its tangent vector.
The second Lagrangian function is just the energy function
\beq
  L_2\left(r,\theta,\frac{dr}{d\lambda}, \frac{d\theta}{d\lambda}\right)
  = \frac12 \left(\frac{dr}{d\lambda}\right)^2+ \frac12 R(r)^2\left(\frac{d\theta}{d\lambda}\right)^2
= E
\eeq
while the first Lagrangian $L_1=ds/d\lambda$ is the speed function
\beq
  L_1\left(r,\theta,\frac{dr}{d\lambda}, \frac{d\theta}{d\lambda}\right)
  = \sqrt{ \left( \frac{dr}{d\lambda} \right)^2+ R(r)^2\left( \frac{d\theta}{d\lambda} \right)^2} 
\,. 
\eeq
Both are independent of the azimuthal angle because of the rotational invariance of the problem
\beq
   R=R(r) \rightarrow \frac{\partial L_i}{\partial\theta} =0\,.
\eeq

The Lagrangian equations are a consequence of extremizing the action
\beq
 0=  \frac{d}{d\lambda} \left(\frac{\partial L_i}{\partial(dx^j/d\lambda)}\right)  
    -\frac{\partial L_i}{\partial x^j} \,.
\eeq
For the second Lagrangian this produces the geodesic equations used in the main body of the article, with the angular equation directly giving the constancy of the angular momentum 
$\ell = \partial L_2 / \partial(d\theta/d\lambda)$. Using instead the first Lagrangian, one obtains the constancy of the momentum conjugate to $\theta$
\begin{eqnarray}\label{eq:ptheta}
  p_\theta &=& \frac{\partial L_1}{\partial(d\theta/d\lambda)}
            = \frac{R^2 d\theta/d\lambda}{\sqrt{(dr/d\lambda)^2 +R^2 (d\theta/d\lambda)^2}}
\nonumber\\
           &=& \frac{R^2 d\theta/d\lambda}{ds/d\lambda}
           = \frac{\ell}{\sqrt{2E}}
 \,,
\end{eqnarray}
where the last line holds not only for an affine parametrization but for any parametrization since the expression is independent of the parametrization, and the final ratio of constants is just the constant vertical component $R(r) \sin\beta$ of the unit speed angular momentum by Eq.~(\ref{eq:unitspeedangmom}).
For the bound geodesics where $|\sin\beta|=1$ at the radial turning point $r=r_{\rm(ext)}$, its absolute value equals the azimuthal radius at the extremal radius: $|p_\theta|=R(r_{\rm(ext)})$.

We can invert the first line of Eq.~(\ref{eq:ptheta}), solving it for the angular velocity $d\theta/d\lambda$ as a function of the canonical momentum $p_\theta$ and then re-express the Lagrangian function $L_1$  in terms of that momentum rather than the angular velocity
\beq
   \frac{d\theta}{d\lambda}
   = \pm \frac{p_\theta\, dr/d\lambda}{R\sqrt{R^2-p_\theta^2}} \,,\quad
   L_1= \frac{ds}{d\lambda}
   = \frac{R |dr/d\lambda|}{\sqrt{R^2-p_\theta^2}}
\,.
\eeq
Suppose we consider geodesic segments for which $r$ is a monotonically increasing or decreasing function between $r_1$ and $r_2$ so $dr/d\lambda\neq0$, allowing it to be used as a parameter itself along the geodesic, i.e., segments not containing any radial turning points except possibly at the endpoints. 
Then one can integrate the previous equations with $r=\lambda$
\begin{eqnarray}
   \theta_2-\theta_1
        &=& \int_{r_1}^{r_2} \frac{d\theta}{dr} \,dr
        =  \pm \int_{r_1}^{r_2} \frac{p_\theta}{R\sqrt{R^2-p_\theta^2}} \,dr 
      \equiv \pm F(r_1,r_2;p_\theta)\,,
\nonumber\\
  s_{21} &=&  \int_{r_1}^{r_2} \frac{ds}{dr} \,dr
         =  \int_{r_1}^{r_2} \frac{R}{\sqrt{R^2-p_\theta^2}} \,dr
 \equiv  G(r_1,r_2;p_\theta)
\,.
\end{eqnarray}
These integrals can be evaluated exactly in terms of elliptic functions, the first as a long combination of such functions, while the second is a relatively short expression in one such function. The first of these two gives the orbit equation $\theta=\theta(r)$ for a given specified momentum $p_\theta$. This relationship cannot be inverted to give $r=r(\theta)$.
However, it can also be thought of as a very complicated condition on the momentum that in turn determines the initial angle that a geodesic has to be shot from its initial point $(r_1,\theta_1)$ to terminate at the final point  $(r_2,\theta_2)$. 
Given that solution for $p_\theta$ one can then evaluate the arc length of the corresponding geodesic using the second integral formula.
For near enough points on the torus, there will be a unique shortest length geodesic which satisfies the previous condition, but for most bound geodesic segments, $r$ is not monotonic for the shortest length geodesic between two given points, but involves one or two radial turning points.

\begin{figure}[p] 
\vglue-0.5in
\typeout{*** EPS figure geodesic spray}
\begin{center}
\includegraphics[scale=0.65]{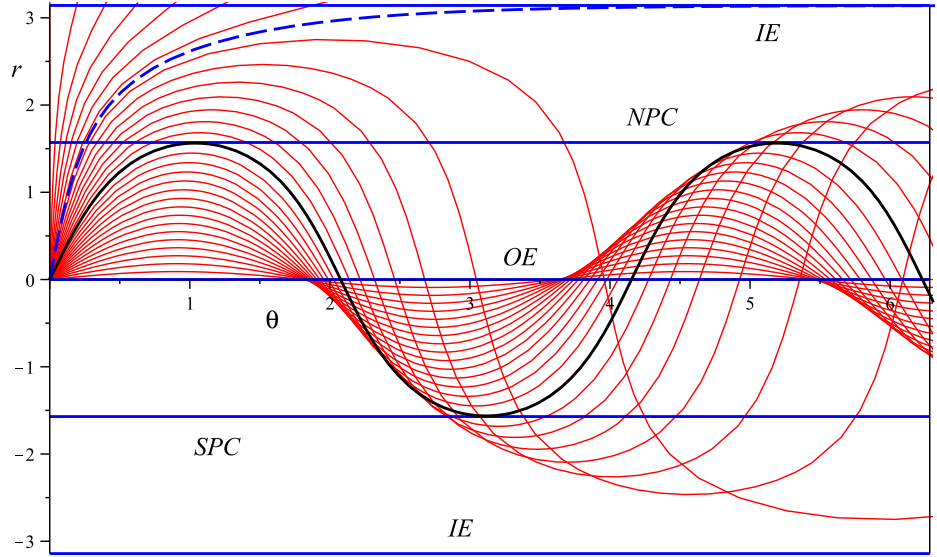}
\end{center}
\vglue-0.1in
\caption{
The spray of geodesics leaving the origin of coordinates in the first quadrant on the unit ring torus
for $-\pi\le r \le\pi,\, 0\le \theta \le 2\pi$.
As one increases the momentum $p_\theta$ from 0, first one passes through the unbound geodesics which cross the inner equator at $r=\pi$ (up to the thick dashed curve which asymptotically approaches the inner equator), then the bound geodesics with a smaller and smaller turning point radius $r_{\rm(max)}$. Looking at where these geodesics cross the Northern Polar Circle at $r=\pi/2$, starting at $\theta=0$ the crossing point moves to the right through unbound geodesics and then into the bound geodesics where their turning point lowers until it reaches that circle (the thick black geodesic, very close to the $[3,2,0]$ closed geodesic which is slightly higher), after which the geodesics which reach that circle overshoot it first, crossing over and then returning to that circle as their second crossing point moves to the right and the turning point rises. Note also the half wavelength $\Delta\theta=\pi/\sqrt{3}\approx 1.81$ of the small oscillations about the outer equator. For $0<\theta<\Delta\theta$ only the outer equator itself from this spray reaches points on that outer equator, but for  $\Delta\theta <\theta<2\Delta\theta$ on the outer equator, a second member of this family reaches the outer equator (with shorter length), while for $2\Delta\theta <\theta<3\Delta\theta$ a third member of this family reaches the outer equator. Of course for $\pi<theta<2\pi$, shorter geodesics arrive from the opposite azimuthal direction.
} 
\label{fig:spray}
\end{figure}

\begin{figure}[t] 
\typeout{*** EPS figure azimuthal angle versus p}
\begin{center}
\includegraphics[scale=0.4]{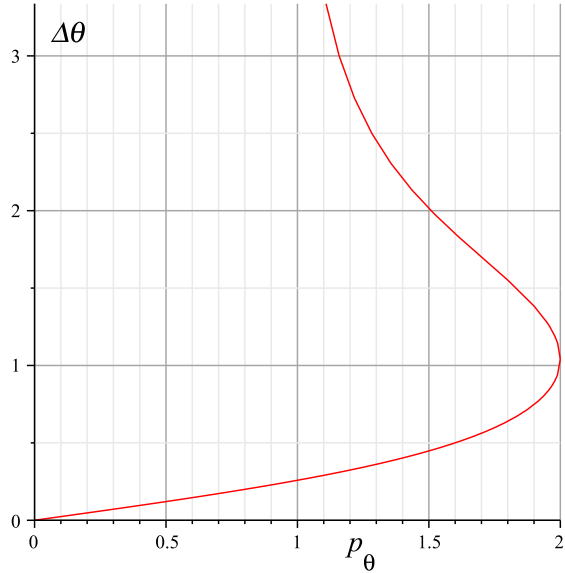}
\end{center}
\vglue-0.3in
\caption{
The plot of $\Delta\theta=\theta_2-\theta_1$ versus the canonical momentum $p_\theta$ for the unit ring torus for geodesics from $r_1=0$ to $r_2=\pi/2$ (outer equator to the Northern Polar Circle), illustrated in the previous figure with $\theta_1=0$. Increasing $p_\theta$ from 0  through the unbound geodesic directions to the bound geodesic directions until the maximal turning point radius $r$ decreases to $\pi/2$,
$\theta_2-\theta_1$ increases to the vertical tangent at $p_\theta=2=R(\pi/2)$ and then the curve turns over as the maximal radius again increases while $p_\theta$ decreases to $1=R(\pi)$ where a vertical asymptote occurs corresponding to the asymptotic approach to the inner equator.
} 
\label{fig:thetaofp}
\end{figure}

This integral function $F$, say when $0\le r_1<r_2<\pi b$ for concreteness, is a monotonically increasing odd function of the momentum $p_\theta$ on the closed interval $[-R(r_2),R(r_2)]$
with value 0 at $p_\theta=0$ and which is concave up when $p_\theta \ge0$, having an extreme value at each endpoint where a vertical tangent occurs when the radial turning point decreases to $r_2$ as $|p_\theta|$ increases to $R(r_2)$. At this point one can decrease $|p_\theta|$ from $R(r_2)$ while allowing the radial turning point to occur within the geodesic segment as that segment extends farther and farther past the original limiting angular interval, corresponding to the curve $\theta_2-\theta_1$ versus $p_\theta$ turning back on itself after the vertical tangent.
Corresponding to this new situation, the geodesic  segment requires overshooting, i.e., going outside the interval from $r_1$ to $r_2$ to reach a radial turning point value $r_{\rm(ext)}$ beyond $r_2$ and return back to it, and the discussion gets much more complicated. See Fig.~\ref{fig:spray}. One must then  break up the azimuthal angle integral into two separate integrals 
\beq
   \theta_2-\theta_1
        = \int_{r_1}^{r_{\rm(max)}} \frac{d\theta}{dr} \,dr
       + \int_{r_2}^{r_{\rm(max)}} \frac{d\theta}{dr} \,dr       \,,
\eeq
where the relation $R(r_{\rm(max)}) = |p_\theta|$ is inverted to give $r_{\rm(max)}$ as a function of $|p_\theta|$
\beq
 R(r_{\rm(max)}) =a+b\cos(r_{\rm(max)}/b) =|p_\theta|\,,\quad
 r_{\rm(max)}=b \arccos\left(\frac{|p_\theta|-a}{b}\right) \,.
\eeq
Now as one decreases $|p_\theta|$ from the value $R(r_2)$, $r_{\rm(max)}$ increases from $r_2$ and $|\theta_2-\theta_1|$ continues increasing as shown in Fig.~\ref{fig:thetaofp} for $p_\theta\ge0$.

This complicates the problem by making the upper limits of integration depend explicitly on the momentum as well, so that it becomes a more complicated condition for that momentum. The two separate integrals correspond respectively to the rise of the radial coordinate to the radial turning point and then its subsequent fall to the final value of the radius. These are also both improper integrals since the denominator goes infinite at the radial turning points, making the problem more tricky numerically as well to find a numerical value of $p_\theta$ which will satisfy the condition. The worst case would involve geodesic segments from the upper hemitorus to the lower hemitorus which require overshooting on both ends of the curve.

\section{Convergence of geodesics}

This is an aside for the experts.
We derived the orbit equation for small geodesic deviations from the outer equator, which can be reparametrized by the arclength $s=(a+b)\,\theta$ along the outer equator
\beq\label{eq:rdeviation}
   r = \pm \mathcal{R} \sin((c+2)^{1/2}\theta) = \pm \mathcal{R} \sin\left(\omega_s \, s\right) \,.
\quad \omega_s = \frac{(c+2)^{1/2}}{a+b}= b^{-1}(c+2)^{-1/2}\,.
\eeq 
These oscillations have the azimuthal angular half-period $T_\theta = \pi/(c+2)^{1/2}$ and the arc length half-period
$L=(a+b)T_\theta = b(c+2)^{1/2}\pi=\pi/\omega_s$.  Note that $L>b\pi$ (the half-circumference of a meridian) and $L<(a+b)\pi=b(c+2) \pi$ (half-circumference of the outer equator)  
for all allowed values of $c\ge -1$ except $c=-1$ which corresponds to the limiting case of a sphere where all three of these lengths coincide. The limit $c\to\infty$, for which $L\to\infty$, corresponds instead to opening up the torus to an infinite flat cylinder, 

This result is just a special case of the geodesic deviation equation \cite{geodeviation} evaluated for deviations from the outer equator geodesic. As illustrated in Fig.~\ref{fig:thetaofp}, we have seen that points on the outer equator separated by a distance greater than $L$  can be connected by a shorter geodesic which deviates from the outer equator. This is an example of the focusing of geodesics by positive curvature. The geodesics which emanate in all directions from an initial point on the outer equator begin crossing each other first at this distance along the outer equator.

One can evaluate the single independent nonzero orthonormal component of the intrinsic curvature tensor (just half the curvature scalar in 2 dimensions \cite{curvaturetensor} and equal to the Gaussian curvature of the surface \cite{gaussiancurvature}) to find \cite{irons}
\beq
   R^{\hat r}{}_{\hat\theta \hat r\hat\theta} = \frac{\cos(r/b)}{b(a+b\cos(r/b))} \,.
\eeq
For ring tori where $a+b\cos(r/b)$ is always positive, this is negative on the inner hemitorus, zero on the Polar Circles and positive on the outer hemitorus,
whose maximum (positive) value occurs at the outer equator
\beq
   R^{\hat r}{}_{\hat\theta \hat r\hat\theta}|_{r=0} = \frac{1}{b(a+b)} =\frac{1}{b^2(c+2)}\,,
\eeq
This is where geodesics cross more quickly than at any other location on the torus, i.e., those which remain close to the outer equator will experience the maximum convergence. This convergence of geodesics by regions of positive curvature leads to the physical application of gravitational lensing.

The unit vector $e_{\hat r}=e_r$ is parallel transported along the outer equator  (as are $e_\theta$ and $e_{\hat\theta}$) while remaining orthogonal to the outer equator. With the geodesic separation vector $\xi=r(s) e_{\hat r}$, the geodesic deviation equation reduces to
\beq
   \frac{d^2 r}{ds} = -R^{\hat r}{}_{\hat\theta \hat r\hat\theta}|_{r=0} \, r
              =  -\omega_s{}^2 \,r\,,
\eeq
which leads to the same solutions (\ref{eq:rdeviation}) that we started with.
Thus we can conclude that for any two points on a nontrivial torus ($c>-1$) separated by less than the radial half-circumference $b\pi$ which in turn is less than the convergence arclength $L$, there will be a unique geodesic of lesser arclength which connects them. Fixing our attention on a given point of the torus, a Euclidean sphere in space of radius $b$ will therefore enclose a region of the surface on which the geodesic from this central point to any other point in this region will be unique.  In other words by evaluating the Euclidean distance between two points on the torus, one has a sufficient condition to guarantee the uniqueness of the boundary value problem for the geodesic equations.
Both leading computer algebra systems Maple and Mathematica have numeric ordinary differential equation solvers which can now handle boundary value problems, so they will produce a geodesic between two specified points. Thus when the two points are within a Euclidean sphere of radius $b$, these programs will deliver the minimal length geodesic.

Finally for spindle tori, the curvature scalar approaches negative infinity as one moves in along a meridian towards to symmetry axis where $a+b\cos(r/b)=0$, and then again turns positive in the lemon part of the surface decreasing from positive infinity to a minimum positive value at the inner equator on the lemon. The curvature scalar there has the positive value
\beq
   R^{\hat r}{}_{\hat\theta \hat r\hat\theta}|_{r=b\pi} = -\frac{1}{b(a-b)} =\frac{1}{b^2(-c)}\,,
\eeq
which is larger than the value at the outer equator. The associated convergence length $b(-c)^{1/2}\pi$
is therefore shorter on the inner equator (but still longer than the half-circumference of the inner equator), but increases as $c$ decreases to the limiting value $-1$ where the inner and outer equators merge into the single equator of the sphere.

\clearpage

\section*{Acknowledgements}

This work would not have been possible without the initial Maple worksheet shared with me by my colleague Klaus Volpert at just the right moment and then enriched by the internet enabled interaction at a distance with Ilteris Murat Derici,
whose interest in the minimal geodesics on the horn torus helped me go further than I had anticipated.

\end{document}